\documentclass[12pt]{amsart}
\usepackage[hmargin=0.8in,height=8.6in]{geometry}
\usepackage{amssymb,amsthm, times}
\usepackage{delarray,verbatim}
\usepackage{srcltx}
\usepackage{ifpdf}
\ifpdf
\usepackage[dvips]{graphicx}
 \else
\usepackage[dvips]{graphicx}
 \fi

\usepackage{ifpdf}
\usepackage{color}
\definecolor{webgreen}{rgb}{0,.5,0}
\definecolor{webbrown}{rgb}{.8,0,0}
\definecolor{emphcolor}{rgb}{0.95,0.95,0.95}

\usepackage{hyperref}
\hypersetup{%
          colorlinks=true,
          linkcolor=webbrown,
          filecolor=webbrown,
          citecolor=webgreen,
          breaklinks=true}
\ifpdf \hypersetup{pdftex,
            pdfstartview=FitH, %%Fit, FitB, FitH
            bookmarksopen=true,
            bookmarksnumbered=true
} \else \hypersetup{dvips} \fi

\numberwithin{equation}{section}

\newtheorem{proposition}{Proposition}[section]

\newtheorem{remark}{Remark}[section]
\newtheorem{lemma}{Lemma}[section]

\newtheorem{assump}{Assumption}[section]

\numberwithin{remark}{section} \numberwithin{proposition}{section}
\numberwithin{corollary}{section}

\renewcommand{\S}{\mathcal{S}}

\newcommand {\ME}{\mathbb{E}^{x}}

\newcommand {\R}{\mathbb{R}}

\newcommand {\F}{\mathcal{F}}
\newcommand {\A}{\mathcal{A}}

\newcommand {\p}{\mathbb{P}}

\newcommand {\E}{\mathbb{E}}

\newcommand{\diff}{{\rm d}}

\newcommand{\word}{\hspace{0.2cm}}
\newcommand{\conn}{\quad\text{and}\quad}
\newcommand{\1}{\mbox{1}\hspace{-0.25em}\mbox{l}}
\newcommand{\lev}{L\'{e}vy }

\newcommand{\s}{s'}
\newcommand{\II}{\mathcal{I}}

\title[Optimal Stopping of Diffusion and its Maximum ]{Explicit Solutions for Optimal Stopping of Linear Diffusion \\and Its Maximum}
\author[M. Egami]{Masahiko Egami}
\address[M. Egami]{Graduate School of Economics,
Kyoto University, Sakyo-Ku, Kyoto, 606-8501, Japan}
\email{egami@econ.kyoto-u.ac.jp}
\urladdr{http://www.econ.kyoto-u.ac.jp/{\textasciitilde}egami/}
\thanks{First Draft: September 28, 2015.  This version: March 31, 2016. This work is in part supported
by Grant-in-Aid for Scientific Research (B) No. 26285069, Japan Society for the Promotion of Science. %An early version of this article was circulated under the title ``Explicit Solutions for Optimal Stopping of Maximum Process %with Absorbing Boundary that Varies with It (https://arxiv.org/pdf/1509.08203v2).
}
\author[T. Oryu]{Tadao Oryu}
\address[T. Oryu]{Institute of Economic Research,
Kyoto University, Sakyo-Ku, Kyoto, 606-8501, Japan}
\email{oryu@kier.kyoto-u.ac.jp}

\date{}
\begin{document}
\begin{abstract}
We provide, in a general setting, explicit solutions for optimal stopping problems that involve  diffusion process and its running maximum.
Our approach is to use the excursion theory for \lev processes. Since general diffusions are, in particular, not of independent increments,  we use an appropriate measure change to make the process have that property. Then we rewrite the original two-dimensional problem as an infinite number of one-dimensional ones and complete the solution. We show general solution methods with explicit value functions and corresponding optimal strategies, illustrating them by some examples.
\end{abstract}

%\1_{\{\tau<\zeta\}}
\maketitle \noindent \small{\textbf{Key words:} Optimal stopping; excursion theory;
 diffusions.\\
%\noindent JEL Classification: G32, D81, C61 \\
 }\\

\section{Introduction}
%We study the following optimal stopping problems.
We let  $X=(X_t, t\ge 0)$ be one-dimensional diffusion and denote by $Y$ the reflected process,
  \[
  Y_t=S_t -X_t
  \]
  where $S_t=\sup _{u \in [0,t]} X_u\vee s$ with $s=S_0$. Hence $Y$ is the excursion of $X$ from its running maximum $S$.  We consider an optimal stopping problem that involves both $X$ and $S$.  It is subject to absorbing boundary that varies with $S$.  That is,
  \begin{eqnarray}\label{eq:main}
\bar{V}(x,s)&=&\sup_{\tau} \E^{x,s} \left[ \int^{\tau}_0 e^{-qt}f(X_t,S_t)\diff t+ e^{-q\tau}g(X_{\tau},S_{\tau}) \right].
\end{eqnarray}
%subject to  absorption \[\zeta:=\inf \{t\geq 0 : S_t-X_t > b(S_t)\}.\]
where the rewards $f$ and $g$ are measurable functions from $\R^2$ to $\R_+$ and $b: \R\mapsto \R_+$ is also a measurable function. The rigorous mathematical definition of this problem is presented in Section \ref{sec:model}.    %This setup means that while $X$ grows and keeps attaining new maxima, the absorbing boundary is accompanying with $S$.
In this paper, we shall solve for optimal strategy and corresponding value function along with optimal stopping region in the $(s, x)$-plane. %The existence of the absorbing state makes the stopping region more complex.

For optimal stopping that involve both $S$ and $X$, we mention a pioneering work of Peskir \cite{peskir1998} where the author used the ``maximum principle".   There are also  Ott \cite{ott_2013} and Guo and Zervos \cite{Guo-Zervos_2010}.  In the former paper, the author solves problems including a capped version of the Shepp-Shiryaev problem (the Russian option) \cite{shepp-shiryaev-1993}.  The latter makes another contribution that extends \cite{shepp-shiryaev-1993} and the reward function therein is $g(x, s)=(x^a s^b-K)^+$ with $a, b, K\ge 0$. This reward function includes perpetual call, lookback option, etc. as special cases.  In many solved problems, Brownian motion or geometric Browninan motion is used as underlying process in an effort to obtain tractable solutions.  A recent development in this area includes Alvarez and Mato\"{a}ki \cite{alvarez-matomaki2014} where a discretized approach is taken to find optimal solutions and a corresponding numerical algorithm is presented.

The idea of our solution method is the following: we look at excursions that occur from each level of $S$, during an excursion from level $S_t=s$, the value of $S_t$ is fixed until $X$ returns to $s$.  Using this, the problems reduce  to an infinite number of one-dimensional optimal stopping problems.  In this spirit, we then attempt to rewrite the problem equation \eqref{eq:main} in the form of sequences of excursions.  For finding the explicit form of the value function, we employ the theory of excursion of \lev processes, in particular the characteristic measure that is related to the height of excursions.  (Refer to Bertoin \cite{Bertoin_1996} as a general reference.)
 Since the diffusion $X$ is not in general of independent increments, we use the measure change \eqref{eq:change-of-measure} to make the diffusion behave like a Brownian motion under the new measure.  Having done that, we solve, at each level of $S$, one-dimensional optimal stopping problems by using the excessive characterization of the value function.  This corresponds to the concavity of the value function after certain transformation. We briefly review this transformation in Section \ref{sec:diffusion-facts}. See Dynkin\cite{dynkin}, Alvarz\cite{alvarez2} and Dayanik and Karatzas \cite{DK2003} for more details.
For the excursion theory for spectrally negative \lev processes (that have only downward jumps), we mention Avram et al. \cite{avram-et-al-2004}, Pistorius \cite{Pistorius_2004} \cite{Pistorius_2007} and Doney \cite{Doney_2005} where, among others, an exit problem of the reflected process $Y$ is  studied.

Our contributions in this paper may advance the literature in several respects:  we do not assume any specific forms or properties in the reward functions (except for mild ones), and we provide \emph{explicit} forms of the value function % \emph{with or without the absorbing boundary}
and illustrate the procedure of the solution method.  In contrast to the literature, our approach is rather direct,  without resorting to the variational inequality.  See Section \ref{sec:solution}.   Accordingly, the solution method can be applied in a more general setting.

The rest of the paper is organized as follows.  In Section
\ref{sec:model}, we formulate a mathematical model with a review of some important facts of linear diffusions, and then
find an optimal solution.  The key step is to represent and compute the value $V(s, s)$, which is handled in Sections \ref{sec:solution} and \ref{sec:finding-V(s,s)}. Under the mild assumptions (Assumption \ref{assumption}), we present $V(s, s)$ in an explicit form.  The next step is to find $V(x, s)$ in Section \ref{sec:method}.   Moreover, we shall demonstrate the methodology by using a new problem (Section \ref{sec:call-with-s}) as well as some problems in the literature (Section \ref{sec:lookback} and \ref{sec:perpetual}).  Let us stress that the new problem might not be easily handled by the conventional methods.  Appendices include a proof of the technical lemma and some discussion about the \emph{smooth-fit principle} in this problem.

\section{Mathematical Model}\label{sec:model} %%%%%%%%%%%%%%%%%%%%%%%%%%%%%%%%%%%%%%%

Let the diffusion process $X=\{X_t;t\geq 0\}$ represent the state variable defined on the probability space
$(\Omega, \F, \p)$, where $\Omega$ is the set of all possible realizations of the
stochastic economy, and $\p$ is a probability measure defined on $\F$. The state space of $X$ is given by
$(l, r):=\mathcal{I}\subseteq \mathbb{R}$, where $l$ and $r$ are \emph{natural boundaries}. That is, $X$ cannot start from and exit from $l$ or $r$.
 We denote by
$\mathbb{F}=\{\F_t\}_{t\ge 0}$ the filtration with respect to which $X$ is adapted and with the usual
conditions being satisfied. We assume that $X$ satisfies the following stochastic differential equation:
\[
\diff X_t=\mu (X_t)\diff t + \sigma (X_t) \diff B_t,\quad X_0=x,
\]
where $\mu (x), \sigma(x) \in \R$ for any $x\in \R$ and $B=\{B_t: t\ge 0\}$ is a standard Brownian motion.

The running maximum process $S=\{S_t;t\geq 0\}$ with $s=S_0$ is defined by
$S_t=\sup _{u \in [0,t]} X_u\vee s$. In addition, we write $Y$ for the reflected process defined by $Y_t=S_t-X_t$.
%and let $\zeta$ be the stopping time defined by
%\[\zeta:=\inf \{t\geq 0 : S_t-X_t > b(S_t)\},\]
%which is the time of absorption.  Note that $b: \R\mapsto \R_+$ is a measurable function and  the level at which the process $X$ is absorbed depends on $S$.
We consider the following optimal stopping problem and the value function $\bar{V}:\R^2 \mapsto \R$ associated with initial values $X_0=x$ and $S_0=s$;
\begin{eqnarray}\label{problem}
\bar{V}(x,s)&=&\sup_{\tau\in\S} \E^{x,s} \left[ \int^{\tau}_0 e^{-qt}f(X_t,S_t)\diff t+ e^{-q\tau}g(X_{\tau},S_{\tau}) \right]
\end{eqnarray}
where $\mathbb{P}^{x,s}(\,\cdot\,):=\mathbb{P}(\,\cdot\,|\,X_0=x, S_0=s)$ and $\E^{x,s}$ is the expectation operator corresponding to $\mathbb{P}^{x,s}$, $q\geq 0$ is the constant discount rate and $\S$ is the set of all $\mathbb{F}$-adapted stopping times.  The payoff is composed of two parts; the running income to be received continuously until stopped, and the terminal reward part.  The running income function
$f:\R^2 \mapsto \R$ is a measurable function that satisfies
\[
\E^{x,s}\left[\int_0^\infty e^{-qt}|f(X_t,S_t)|\diff t\right]<\infty.
\]
The reward function $g:\R^2 \mapsto \R_+$ is assumed to be measurable. Our main purpose is to calculate $\bar{V}$ and to find the stopping time $\tau^*$ which attains the supremum.

For each collection $D=(D(s))_{s\in\R}$ of Borel measurable sets, we define a stopping time $\tau(D)$ by
\begin{equation}\label{eq:tau-D}
\tau(D):=\inf\{t\geq 0:S_t-X_t\in D(S_t)\},
\end{equation}
and define a set of stopping times $\S'$ by $\S':=\{\tau(D)\}$. In other words, $\tau(D)$ is the first time the excursion $S-X$ from level, say $S=s$, enters the region $D(s)$. Suppose that $D(m) = (c, +\infty)$ for any $m\in\R$, we write
\[
\tau_c:=\inf\{t\geq 0:S_t-X_t>c\}.
\]

\subsection{Optimal Strategy}\label{subsec:optimality}  %%%%%%%%%%%%%%%%%%%%%%%%%%%%%%%%%%%%%%%%
We will reduce the problem \eqref{problem} to an infinite number of one-dimensional optimal stopping problem and discuss the optimality of the
proposed strategy \eqref{eq:tau-D}.  Let us denote by $\bar{f}:\R^2\mapsto\R$ the $q$-potential of $f$:
\[
\bar{f}(x,s):=\E^{x,s} \left[ \int^{\infty}_0 e^{-qt}f(X_t,S_t)\diff t \right].
\]
From the strong Markov property of $(X,S)$, we have
\begin{eqnarray}\label{eq:potential-rewrite}
&&\E^{x,s} \left[ \int^{\tau}_0 e^{-qt}f(X_t,S_t)\diff t \right]\\
&=&\E^{x,s} \left[\int^{\infty}_0 e^{-qt}f(X_t,S_t)\diff t - \int^{\infty}_{\tau} e^{-qt}f(X_t,S_t)\diff t \right]\nonumber\\
&=&\bar{f}(x,s)-\E^{x,s} \left[\E \left[\int^{\infty}_{\tau} e^{-qt}f(X_t,S_t)\diff t \Bigm| \mathcal{F}_{\tau}\right] \right]\nonumber\\
&=&\bar{f}(x,s)-\E^{x,s} \left[e^{-q\tau}\E^{X_{\tau},S_{\tau}}\left[\int^{\infty}_0 e^{-qt}f(X_t,S_t)\diff t \right] \right]\nonumber\\
&=&\bar{f}(x,s)-\E^{x,s} \left[e^{-q\tau}\bar{f}(X_{\tau},S_{\tau})\right]\nonumber
%&=&\bar{f}(x,s)-\E^{x,s} \left[e^{-q\tau}\bar{f}(X_{\tau},S_{\tau})\right].\nonumber
\end{eqnarray}
Hence the value function $\bar{V}$ can be written as
\begin{eqnarray}\label{eq:value-function-rewritten}
\bar{V}(x,s)&=&\bar{f}(x,s)+V(x,s),\nonumber
\end{eqnarray}
where
\begin{align}\label{eq:V-modified}
\hspace{0.5cm} V(x,s):=\sup_{\tau\in\S}\E^{x,s}\left[e^{-q\tau}(g-\bar{f})(X_{\tau},S_{\tau})\right].
\end{align}
Since $\bar{f}(x, s)$ has nothing to do with the choice of $\tau$, we concentrate on $V(x, s)$.

Let us first define the first passage times of $X$:
\[T_a:=\inf\{t\geq 0:X_t > a\}\quad\text{and}\quad T_a^{-}:=\inf\{t\ge 0: X_t
 <a\}.\]
By the dynamic programming principle, we can write $V(x,s)$ as
\begin{align}\label{eq:dp}
V(x,s)=\sup_{\tau\in\S}\E^{x,s}\left[\1_{\{\tau<\theta\}}e^{-q\tau}(g-\bar{f})(X_{\tau},S_{\tau})+\1_{\{\theta<\tau\}}e^{-q\theta}V(X_{\theta},S_{\theta})\right],
\end{align}
for any stopping time $\theta\in\S$. See, for example, Pham \cite{Pham-book} page 97. Now we set $\theta =T_s$ in (\ref{eq:dp}).
For each level $S=s$ from which an excursion $Y=S-X$ occurs, the value $S$ does not change during the excursion.
Hence, during the first excursion interval  from $S_0=s$, we have $S_t=s$ for any $t\leq T_s$, and (\ref{eq:dp}) can be written as the following one-dimensional problem for the state process $X$;
\begin{eqnarray}\label{eq:one-dim-version}
V(x,s)&=&\sup_{\tau\in\S}\E^{x,s}\left[\1_{\{\tau<T_s\}}e^{-q\tau}(g-\bar{f})(X_{\tau},s)
+\1_{\{T_s<\tau\}}e^{-qT_s}V(s,s)\right].\nonumber
\end{eqnarray}
Now we can look at \emph{only} the process $X$ and find $\tau^*\in \S$.

In relation to \eqref{eq:one-dim-version}, we consider the following one-dimensional optimal stopping problem as for $X$ and its value function $\widehat{V}:\R^2\mapsto \R$;
\begin{align}\label{eq:Vhat}
\widehat{V}(x,s)&=\sup_{\tau\in\S}\E^{x,s}\left[\1_{\{\tau<T_s\}}e^{-q\tau}(g-\bar{f})(X_{\tau},s)
+\1_{\{T_s<\tau\}}e^{-qT_s}K\right],
\end{align}
where $K\geq 0$ is a constant. Note that $V=\widehat{V}$ holds when $K=V(s,s)$, and we shall present how to characterize and compute $V(s,s)$ in Sections \ref{sec:solution} and \ref{sec:finding-V(s,s)}.

Recall that, in pursuit of optimal strategy $\tau^*$ in the linear diffusion case, we can utilize the full characterization of the value function and  optimal stopping
rule.  In particular, an optimal stopping is,  in a very general setup, attained by the \emph{threshold strategy}.  That is, one should stop $X$ when at the first time that $X$ enters a certain region, called stopping region.
See Dayanik and Karatzas \cite{DK2003}; Propositions 5.7 and 5.14. See also Pham \cite{Pham-book}; Section 5.2.3.  Note that, however, writing the value of $V(s, s)$ in an explicit form is not trivial and is an essential part of the solution, which we shall do in the next section (see Propositions \ref{prop:1} and \ref{prop:2}).

\subsection{Important Facts of Diffusions}\label{sec:diffusion-facts}
Let us recall the fundamental facts about one-dimensional diffusions; let the differential operator $\A$ be the infinitesimal generator of the process $X$ defined by
\[
\A v(\cdot)=\frac{1}{2}\sigma^2(\cdot)\frac{\diff^2  v}{\diff x^2}(\cdot)+\mu(\cdot)\frac{\diff v}{\diff x}(\cdot)
\]
and consider the ODE $(\A-q)v(x)=0$. This equation has two
fundamental solutions: $\psi(\cdot)$ and $\varphi(\cdot)$.
We set
$\psi(\cdot)$ to be the increasing and $\varphi(\cdot)$ to be the
decreasing solution.  They are linearly independent positive
solutions and uniquely determined up to multiplication. It is well
known that
\begin{align}\label{eq:psi-phi}
  \ME[e^{-\alpha\tau_z}]%1_{\{\tau_z<\infty\}}
  =\begin{cases}
\frac{\psi(x)}{\psi(z)}, & x \le z,\\[4pt]
\frac{\varphi(x)}{\varphi(z)}, &x \ge z.
  \end{cases}
\end{align} For the complete characterization of $\psi(\cdot)$ and
$\varphi(\cdot)$, refer to It\^{o} and McKean \cite{IM1974}. Let us now define
\begin{align} \label{eq:F}
F(x)&:=\frac{\psi(x)}{\varphi(x)}, \hspace{0.5cm} x\in
\mathcal{I}.
\end{align}
Then $F(\cdot)$ is continuous and strictly increasing.  Next,
following Dynkin (pp.\ 238, \cite{dynkin}), we define concavity
of a function with respect $F$ as follows:
%Let $F :\mathcal{I}=(c, d)\rightarrow\mathbb{R}$ be a strictly
%increasing function.
A real-valued function $u$ is called \emph{$F$-concave} on $\mathcal{I}$
if, for every $x\in[l, r]\subseteq \mathcal{I}$,
\[
u(x)\geq
u(l)\frac{F(r)-F(x)}{F(r)-F(l)}+u(r)\frac{F(x)-F(l)}{F(r)-F(l)}.\]

Now consider the optimal stopping problem:
\[V(x)=\sup_{\tau\in \S}\E^x[e^{-q \tau}h(X_\tau)]
\]
where $h$: $[c, d]\mapsto \R_+$. Let $W(\cdot)$ be the smallest nonnegative concave majorant of
\begin{equation}\label{eq:transform}
H:=\frac{h}{\varphi}\circ F^{-1} \quad \word \text{on}\word [F(c), F(d)]
\end{equation}
where $F^{-1}$ is the inverse of $F$. Then we have $V(x)=\varphi(x)W(F(x))$ and the optimal stopping region $\Gamma$ is
\[
\Gamma:=\{x\in [c, d]: V(x)=h(x)\} \conn \tau^*:=\inf\{t\ge 0: X_t\in \Gamma\}
\] as in Propositions 4.3 and 4.4 of \cite{DK2003}.  Note that for the rest of this article, the term ``transformation" should be understood as \eqref{eq:transform}.

When both boundaries $l$ and $r$ are natural, $V(x)<+\infty$ for all $x\in (l, r)$ if and only if
\begin{align}\label{eq:finiteness}
  \xi_l:=\limsup_{x\downarrow l}\frac{h^+(x)}{\varphi(x)} \conn \xi_r:=\limsup_{x\uparrow r}\frac{h^+(x)}{\psi(x)}
\end{align}
are both finite.
\section{Representation of $V(s, s)$}\label{sec:solution}
Now we look to an explicit solution of $\bar{V}$ for $\tau\in\S'$. The first step is to find $V(s, s)$ in \eqref{eq:one-dim-version}.
%\subsection{When \bf{$S_0=X_0$}}\label{sec:main}
That it, we consider the case $S_0=X_0$. Set stopping times
$T_m$ as $T_m=\inf\{t\geq 0:X_t > m\}$ and define the function $l_D:\R_+\mapsto\R_+$ by
\begin{equation}\label{eq:lD}
l_D(m):=\inf D(m).
\end{equation}
 %\quad D\in S'
for which $\tau(D) \in \S'$.   Since optimal strategy belongs to threshold strategies, we shall focus on the set of $\tau(D)$ in \eqref{eq:tau-D}.  In other words, if $S_0=X_0$, given a threshold strategy $\tau=\tau(D)$ where  $D(m)$ is in the form of $[a, c]\subset [0, +\infty)$, the value $l_{D(m)}$ is equal to $a$.  Accordingly,
$S_{\tau_{l_D(m)}}=S_{\tau_a}$ on the set $\{S_{\tau_{l_D(m)}}\in\diff m\}$.
%S_{\tau_{l_D(m)}}-X_{\tau_{l_D(m)}}\leq b(m),

From the strong Markov property of $(X,S)$, when $\tau(D)\in\S'$ and
$S_0=X_0=s$, we have
\begin{eqnarray}\label{eq:interim}
&&\E^{s,s}\left[\1_{\{\tau(D)<+\infty\}}e^{-q\tau(D)}(g-\bar{f})(X_{\tau(D)},S_{\tau(D)})\right]\\
&=&\int^{\infty}_s\E^{s,s}\left[\1_{\{\tau(D)<+\infty, S_{\tau(D)}\in\diff m\}}e^{-q\tau(D)}(g-\bar{f})(X_{\tau(D)},S_{\tau(D)})\right]\nonumber\\
&=&\int^{\infty}_s\E^{s,s}\left[\1_{\{T_m\leq\tau(D)\}}e^{-qT_m}\E^{m,m}\left[e^{-q\tau_{l_D(m)}}(g-\bar{f})(X_{\tau_{l_D(m)}},S_{\tau_{l_D(m)}})\right.\right.\nonumber\\
&&\left.\left.\times\1_{\{ S_{\tau_{l_D(m)}}\in \diff m\}}\right]\right]\nonumber\\
&=&\int^{\infty}_s\E^{s,s}\left[\1_{\{S_{\tau(D)}\geq m\}}e^{-qT_m}\right](g-\bar{f})(m-l_D(m),m)\nonumber\\
&&\times\E^{m,m}\left[e^{-q\tau_{l_D(m)}}\1_{\{S_{\tau_{l_D(m)}}\in \diff m\}}\right]\nonumber.
\end{eqnarray}
Now we calculate these expectations by changing probability measure.
We introduce the probability measure $\p^{\varphi}_{x,s}$ defined by
\begin{equation}\label{eq:change-of-measure}
\p^{\varphi}_{x,s}(A):=\frac{1}{\varphi(x)}\E^{x,s}\left[ e^{-qt}\varphi(X_t)\1_{A} \right], \text{for every } A\in\F.
\end{equation}
Then $F(X)=(F(X_t))_{t\in \R_+}$ is a process in natural scale under this measure $\p^{\varphi}_{x,s}$. See Borodin and Salminen \cite{borodina-salminen}(pp.\ 33) and Dayanik and Karatzas \cite{DK2003} for detailed explanations. Hence we can write, under $\p^{\varphi}_{x,s}$, $F(X_t)=\sigma_F B^{\varphi}_t$, where $\sigma_F>0$ is a constant and $B^{\varphi}$ is a Brownian motion under $\p^{\varphi}_{x,s}$. Since $F(X)$ is a \lev process, we can define the process $\eta:=\{\eta_t;t\geq 0\}$ of the height of the excursion as
\[
\eta_u:= \sup\{(S-X)_{T_{u-}+w} : 0\leq w \leq T_u - T_{u-}\}, \text{ if \;}  T_u > T_{u-},
\]
and $\eta_u=0$ otherwise, where $T_{u-}:=\inf\{t\ge0:X_t\ge u\}=\lim_{m\rightarrow u-}T_{m}$.
Then $\eta$ is a Poisson point process, and we denote its
characteristic measure under $\p^{\varphi}_{x,s}$ by $\nu:\F\mapsto\R_+$ of $F(X)$.
It is well known that
\[
\nu[u,\infty)=\frac{1}{u}, \text{ for } u\in\R_+.
\] See, for example, \c{C}inlar \cite{cinlar} (pp.\ 416).
By using these notations, we have\footnote{Note that when the diffusion $X$ is a standard Brownian motion $B$, then $F(x)=x$ and the right-hand side reduces to $\exp\left(-\int_s^m\frac{\diff u}{l_D(u)}\right)$.}
\begin{eqnarray}\label{eq:probability}
\p^{\varphi}_{s,s}(S_{\tau(D)}>m)&=&\exp\left( -\int^{F(m)}_{F(s)}\nu[y-F(F^{-1}(y)-l_D(F^{-1}(y))),\infty)\diff y \right)\nonumber\\
&=&\exp\left( -\int^m_{s}\frac{F'(u)\diff u}{F(u)-F(u-l_D(u))} \right).
\end{eqnarray}
On the other hand, from the definition of the measure $\p^{\varphi}_{x,s}$, we have
\begin{eqnarray*}
\p^{\varphi}_{s,s}(S_{\tau(D)}>m)
&=&\frac{1}{\varphi(s)}\E^{s,s}\left[ e^{-qT_m}\varphi(X_{T_m})\1_{\{ S_{\tau(D)}>m \}} \right]\\
&=&\frac{\varphi(m)}{\varphi(s)}\E^{s,s}\left[ e^{-qT_m}\1_{\{ S_{\tau(D)}>m \}} \right].
\end{eqnarray*}
Combining these two things together,
\begin{equation}
\E^{s,s}\left[ e^{-qT_m}\1_{\{ S_{\tau(D)}>m \}} \right]=\frac{\varphi(s)}{\varphi(m)}\exp\left( -\int^m_s\frac{F'(u)\diff u}{F(u)-F(u-l_D(u))} \right).
\end{equation}
Similarly, by changing the measure and noting that $X_{\tau_{l_D(m)}}=m-l_D(m)$, we have
\begin{eqnarray}
&&\E^{m,m}\left[e^{-q\tau_{l_D(m)}}\1_{\{S_{\tau_{l_D(m)}}\in \diff m\}}\right]\\ \nonumber
&=&\frac{\varphi(m)}{\varphi(m-l_D(m))}\cdot\frac{1}{\varphi(m)}\E^{m,m}\left[e^{-q\tau_{l_D(m)}}\varphi(X_{\tau_{l_D(m)}})\1_{\{S_{\tau_{l_D(m)}}\in \diff m\}}\right]\\ \nonumber
&=&\frac{\varphi(m)}{\varphi(m-l_D(m))}\p^{\varphi}_{m,m}(F(S_{\tau_{l_D(m)}})\in\diff F(m))\\ \nonumber
&=&\frac{\varphi(m)}{\varphi(m-l_D(m))}\p^{\varphi}_{m,m}(F(S_{\tau_{l_D(m)}})-F(X_{\tau_{l_D(m)}})=l_D(m), F(S_{\tau_{l_D(m)}})\in\diff F(m))\nonumber.
\end{eqnarray}

Since $F(X)$ is a \lev process under  $\p^{\varphi}_{s,s}$, we can apply Theorem 2 in Pistorius \cite{Pistorius_2007} to calculate the last probability. Then we have
\begin{equation}\label{eq:reach-m}
\p^{\varphi}_{m,m}(F(S_{\tau_{l_D(m)}})-F(X_{\tau_{l_D(m)}})=l_D(m), F(S_{\tau_{l_D(m)}})\in\diff F(m))=\frac{F'(m)\diff m}{F(m)-F(m-l_D(m))}.
\end{equation}

We have, up to this point, proved the following:
\begin{proposition}\label{prop:1}
When $S_0=X_0$, the function $V(s,s)$ for $\tau\in \S$ can be represented by
\begin{eqnarray}\label{eq:V(s,s) integral form}
V(s,s)&=&\sup_{l_D}\int^\infty_s\frac{\varphi(s)}{\varphi(m-l_D(m))}\exp\left( -\int^m_s\frac{F'(u)\diff u}{F(u)-F(u-l_D(u))} \right)\\
&&\times\frac{F'(m)(g-\bar{f})(m-l_D(m),m)}{F(m)-F(m-l_D(m))}\diff m. \nonumber
\end{eqnarray}
\end{proposition}
This representation of $V(s, s)$ applies to general cases.

\section{Computing $V(s, s)$}\label{sec:finding-V(s,s)}
In solving an optimal stopping problem involving $S$ and $X$, one of the aims is to draw a diagram like Figure \ref{Fig-example-design}.  Note that we draw the diagram with $s$ on the horizontal axis since it is better understood than otherwise.  For distinct values in the $(s, x)$-diagram, we need to determine whether a point in $\R^2$ is in the continuation region ($\mathrm{C}$) or stopping region ($\Gamma$).
 The task in this subsection\footnote{Once this is done, then the next task is to examine the points $(s, x)$ by moving downwards to $x=0$ from the diagonal $x=s$. We take this in Section \ref{sec:method} and solve an example in \ref{sec:call-with-s}.} is to compute the value $V(s, s)$ at a point $(s, s)$ on the diagonal and to determine whether it belongs to $\mathrm{C}$ or $\Gamma$.

 As stressed before, once we fix $S=s$, the problem reduces to one-dimensional problems in $X$ and hence the way to find $V(s, s)$ is similar to the one described in Section \ref{sec:diffusion-facts}. Specifically, we set $S=s$ and attempt to find the smallest nonnegative concave majorant $W_s(\cdot)$ of
\begin{equation}\label{eq:H_s}
H_s(y):=\frac{(g-\bar{f})(F^{-1}(y), s)}{\varphi(F^{-1}(y))}, \quad y\in F(\II)\end{equation}
in the neighborhood of $s$.  Let us denote by $\Sigma_s\subseteq \R$ (resp. $\mathcal{C}_s$) the stopping region (resp. continuation region) with respect to the reward $H_s(y)$ in \eqref{eq:H_s}, corresponding to this $s$. Note that this $\Sigma_s$ (resp. $\mathcal{C}_s$) should be distinguished from the stopping region $\Gamma\subseteq \R^2$ (resp. $\mathrm{C}$) of the problem \eqref{problem}, the final object to figure out.
 Due to the dependence of the reward on $s$, however, there are some situations where we need to be careful.  To discuss, we shall hereafter (for the rest of this section) assume the following:
\begin{assump}\label{assumption}
Denote by $x^*(s)$ the threshold point that separates $\mathcal{C}_s$ and $\Gamma_s$ with respect to the reward $H_s$ associated with this $s$.
\begin{enumerate}
  \item [(i)]  $(g-\bar{f})(x, s)$ is increasing in $s$,
  \item [(ii)] for $s\in \II=(l, r)$, the continuation region $\mathcal{C}_s$ corresponding to $H_s(y)$ is in the form of $A_s:=(l, x^*(s)]$ or $B_s:=[x^*(s), r)$, and
  \item [(iii)] \eqref{eq:finiteness} holds.
\end{enumerate}
\end{assump}
The first assumption is merely to restrict our problems  to practical ones because we are solving maximization problems.  The second is to make our argument concrete and simplistic. More complicated structure can be handled by some combinations of the cases prescribed below.  For the third, since our main concern is to find a finite value function, we shall consider the case where \eqref{eq:finiteness} holds.

\subsection{Case (1): \emph{$s \in \Sigma_s$}}
If $s$ belongs to $\Sigma_s$, we have  $W_s(y)=H_s(y)$; however,  instead of stopping immediately, there is a possibility that a greater value can be attained if one stops $X$ during the excursion from some upper level $s'>s$. Recall that Proposition \ref{prop:1} has incorporated this.  See also Remark \ref{rem:meaning} below.  Now we wish to obtain more explicit formulae for $V(s,s)$ from the general representation \eqref{eq:V(s,s) integral form}.
For this purpose, let us denote
\[
P(u; l_D):=\frac{F'(u)}{F(u)-F(u-l_D(u))}, \conn G(u; l_D):=(g-\bar{f})(u-l_D(u),u),
\]
to avoid the long expression and
rewrite \eqref{eq:V(s,s) integral form} in the following way:
for any $m-s>\epsilon > 0$,
\begin{eqnarray*}
V(s, s)&=&\sup_{l_D}\left[\int_{s}^{s+\epsilon}\frac{\varphi(s)}{\varphi(m-l_D(m))}\exp\left(-\int_s^{m}P(u; l_D)\diff u\right)P(m; l_D)G(m; l_D)\diff m\right.\\
&&+\left.\frac{\varphi(s)}{\varphi(s+\epsilon)}\exp\left(-\int_s^{s+\epsilon}P(u; l_D)\diff u\right)\right.\\
&&\hspace{0.8cm}\left.\times\int_{s+\epsilon}^\infty\frac{\varphi(s+\epsilon)}{\varphi(m-l_D(m))}\exp\left(-\int_{s+\epsilon}^{m}P(u; l_D)\diff u\right)P(m; l_D)G(m; l_D)\diff m\right]\\
&=&\sup_{l_D}\left[\int_{s}^{s+\epsilon}\frac{\varphi(s)}{\varphi(m-l_D(m))}\exp\left(-\int_s^{m}P(u; l_D)\diff u\right)P(m; l_D)G(m; l_D)\diff m\right.\\
&&+\left.\frac{\varphi(s)}{\varphi(s+\epsilon)}\exp\left(-\int_s^{s+\epsilon}P(u; l_D)\diff u \right) V(s+\epsilon, s+\epsilon)\right].
\end{eqnarray*}
This expression naturally motivates us to set $V_\epsilon : \R\mapsto\R$ as
\begin{equation}\label{eq:V-epsilon}
V_\epsilon(s):=\sup_{l_D(s)}\left[ \frac{\varphi(s)}{\varphi(s+\epsilon)}\exp\left(-\epsilon P(s; l_D)\right)V(s+\epsilon,s+\epsilon)+\frac{\varphi(s)}{\varphi(s-l_D(s))}\cdot\epsilon P(s; l_D)G(s; l_D)\right]
\end{equation} and we have $\lim_{\epsilon\downarrow 0}V_\epsilon (s)=V(s,s)$.
Dividing both sides by $\varphi(s)$ and choosing the optimal level $l_D^*(s)$,  we have, for any $s$ given,
\begin{equation}\label{eq:V-epsilon-in-F}
  \frac{V_\epsilon(s)}{\varphi(s)}=\left[\frac{V(s+\epsilon, s+\epsilon)}{\varphi(s+\epsilon)}e^{-\epsilon P(s; l^*_D)}+\frac{\epsilon P(s; l^*_D)G(s; l^*_D)}{\varphi(s-l^*_D(s))}\right].
\end{equation}

Let us consider the transformation of a Borel function $z$ defined on $-\infty\le c\le x\le d\le \infty$ through
\begin{equation}\label{eq:F-trans}
Z(y):=\frac{z}{\varphi}\circ F^{-1}(y)
\end{equation} on $[F(c), F(d)]$ where $F^{-1}$ is the inverse of the strictly increasing $F(\cdot)$ in \eqref{eq:F}.
 If we evaluate $Z$ at $y=F(x)$, we obtain $Z(F(x))=\frac{z(x)}{\varphi(x)}$, which is the form that appears in \eqref{eq:V-epsilon-in-F}. Note that
 \begin{equation}\label{eq:derivative}
  Z'(y)=q'(x) \quad \text{where}\word q'(x)=\frac{1}{F'(x)}\left(\frac{z}{\varphi}\right)'(x)
 \end{equation}

%To make explicit calculations possible, we consider the case where the reward increases as $S$ does, a natural problem formulation.
\begin{proposition}\label{prop:2}
Fix $s\in \mathcal{I}$.  If (1) the reward function $(g-\bar{f})(x, s)$ is increasing in the second argument and (2) if $\log \varphi(\cdot)$ is strictly convex on $\II$, we have
  \begin{equation}\label{eq:V-explicit}
V(s,s)=\frac{\varphi(s)}{\varphi(s-l^*_D(s))}\cdot Q(s; l_D^*) \cdot (g-\bar{f})(s-l^*_D(s),s),
\end{equation}
where
\[
Q(u; l_D):=\frac{F'(u)\varphi'(u)}{\varphi''(u)[F(u)-F(u-l^*_D(u))]+F'(u)\varphi'(u)}
\]
and $l^*_D(s)$ is the maximizer of the map \[z \mapsto
\frac{\varphi(s)}{\varphi(s-z)}\cdot\frac{F'(s)\varphi'(s)}{\varphi''(s)[F(s)-F(s-z)]+F'(s)\varphi'(s)}\cdot(g-\bar{f})(s-z,s)\quad \text{on}\word [0, \infty).\]
If (2)' $\log \varphi(\cdot)$ is linear on $\II$, then $Q(\cdot; l_D)$ is to be replaced by
%\begin{equation}\label{eq:V-explicit-Bm}
%V(s,s)=\frac{\varphi(s)}{\varphi(s-l^*_D(s))}\cdot \tilde{Q}(s; l_D^*) \cdot (g-\bar{f})(s-l^*_D(s),s),
%\end{equation}
%where
\[
\tilde{Q}(u; l_D):=\frac{F'(u)\varphi(u)}{(\varphi'(u)-\varphi(u))[F(u)-F(u-l^*_D(u))]+F'(u)\varphi(u)}.
\]

\end{proposition}
\begin{remark}\normalfont
The strict convexity of $\log \varphi(\cdot)$ implies
\begin{equation}\label{eq:log-varphi}
\frac{\varphi(s)}{\varphi(s+\epsilon)}\frac{\varphi'(s+\epsilon)}{\varphi'(s)}<1\quad  \forall s\in \II \conn \forall \epsilon>0.
\end{equation}
On the other hand, if $\log \varphi(\cdot)$ is linear, then the inequality in \eqref{eq:log-varphi} is replaced by the equality:
\begin{equation}\label{eq:log-varphi2}
\frac{\varphi(s)}{\varphi(s+\epsilon)}\frac{\varphi'(s+\epsilon)}{\varphi'(s)}=1\quad  \forall s\in \II \conn \forall \epsilon>0.
\end{equation}
We shall use this property in the  proof of Lemma \ref{lem:epsilon}.  Note that it is easily proved that geometric Brownian motion satisfies \eqref{eq:log-varphi} and Brownian motion (with or without drift) does the other one.  The $\varphi(\cdot)$ functions of Orstein-Uhlenbeck process $\diff X_t=k(m-X_t)\diff t+\sigma \diff B_t$ (with $k>0, \sigma>0$ and $m\in \R$) and its exponential version $\diff X_t=\mu X_t(\alpha-X_t)\diff t+\sigma X_t\diff B_t$ (where $\mu, \alpha, \sigma$ are positive constant) involve special functions; the parabolic cylinder function and the confluent hypergeometric function of the second kind, respectively (see Lebedev \cite{lebedev}).   While it is hard to prove the log convexity of these special functions, it is numerically confirmed that both processes satisfy \eqref{eq:log-varphi}.
Moreover, the equality condition \eqref{eq:log-varphi2} implies that $(\log \varphi(s+\epsilon))'=(\log \varphi(s))'$ so that $\varphi(s)$ is an exponential function.  Hence it is confirmed that this case includes Brownian motion.
\end{remark}

\begin{proof}(of Proposition \ref{prop:2})
Let us first take the strictly convex case of $\log \varphi(\cdot)$.  For taking limits of $\epsilon\downarrow 0$ in \eqref{eq:V-epsilon}, we need the following lemma whose proof is postponed to Appendix \ref{app:A}:
\begin{lemma}\label{lem:epsilon} Under the assumption of Proposition \ref{prop:2} with convex $\log\varphi(\cdot)$, for $\epsilon>0$ sufficiently close to zero, we have
\begin{equation} \label{eq:V-V-1}
    \frac{V_\epsilon(s)}{\varphi(s)}=\alpha_s(\epsilon)\cdot\frac{V(s+\epsilon, s+\epsilon)}{\varphi(s+\epsilon)} \quad \word\text{where}\word
    \alpha_s(\epsilon):=\frac{\varphi'(s+\epsilon)}{\varphi'(s)}.
    %\frac{V_\epsilon(s)}{\varphi(s)}=\left(1+\epsilon\cdot\frac{\varphi''(s)}{\varphi'(s)}\right)\frac{V(s+\epsilon)}{\varphi(s+\epsilon)}
\end{equation}
\end{lemma}
\noindent Note that $\alpha_s(\epsilon)\in (0, 1)$ for all $s\in \mathcal{I}$ and $\epsilon>0$ and that $\alpha_s(\epsilon)\uparrow 1$ for all $s\in \mathcal{I}$.\\

Suppose that the lemma is proved, let us continue the proof of Proposition \ref{prop:2}. By using \eqref{eq:V-V-1} in Lemma \ref{lem:epsilon}, we can write, for $\epsilon$ small,
\begin{equation}\label{eq:V_epsilon}
V_\epsilon(s)-\frac{\varphi(s)}{\varphi(s+\epsilon)}\exp\left(-\epsilon P(s; l^*_D)\right)V(s+\epsilon, s+\epsilon)=\left(1-\frac{\varphi'(s)}{\varphi'(s+\epsilon)}\exp\left(-\epsilon P(s; l^*_D)\right)\right)V_\epsilon(s).
\end{equation}
Moreover, since $\lim_{\epsilon\downarrow 0}V(s+\epsilon,s+\epsilon)=V(s,s)$, the optimal threshold $l^*_D(s)$ should satisfy
\begin{equation*}
V(s, s)=\lim_{\epsilon\downarrow 0}V_\epsilon(s)=\lim_{\epsilon\downarrow 0}\left[ \frac{\varphi(s)}{\varphi(s+\epsilon)}\exp\left(-\epsilon P(s; l^*_D)\right)V(s+\epsilon,s+\epsilon)
+\frac{\varphi(s)}{\varphi(s-l^*_D(s))}\cdot\epsilon P(s; l^*_D)G(s; l^*_D)\right],
\end{equation*}
%\begin{eqnarray*}
%\lim_{\epsilon\downarrow 0}V_\epsilon(s)&=&\lim_{\epsilon\downarrow 0}\left[ \frac{\varphi(s)}{\varphi(s+\epsilon)}\exp\left(-\frac{\epsilon F'(s)\diff u}{F(s)-F(s-l^*_D(s))}\right)V(s+\epsilon,s+\epsilon)\right.\\
%&&\left.+\frac{\varphi(s)}{\varphi(s-l^*_D(s))}\cdot\frac{\epsilon F'(s)(g-\bar{f})(s-l^*_D(s),s)}{F(s)-F(s-l^*_D(s))} \right].
%\end{eqnarray*}
from which equation, in view of \eqref{eq:V_epsilon}, we obtain
\begin{eqnarray*}
V(s,s)&=&\lim_{\epsilon\downarrow 0}\frac{ V_\epsilon(s)-\frac{\varphi(s)}{\varphi(s+\epsilon)}\exp\left(-\epsilon P(s; l^*_D)\right)V(s+\epsilon,s+\epsilon) }
{ 1-\frac{\varphi'(s)}{\varphi'(s+\epsilon)}\exp\left(-\epsilon P(s; l^*_D)\right) }\\
&=&\lim_{\epsilon\downarrow 0}\frac{ \frac{\varphi(s)}{\varphi(s-l^*_D(s))}\cdot\epsilon P(s; l^*_D)G(s; l^*_D) }
{ 1-\frac{\varphi'(s)}{\varphi'(s+\epsilon)}\exp\left(-\epsilon P(s; l^*_D)\right) }\\
&=&\frac{\varphi(s)}{\varphi(s-l^*_D(s))}
\frac{F'(s)\varphi'(s)}{\varphi''(s)[F(s)-F(s-l^*_D(s))]+F'(s)\varphi'(s)}
(g-\bar{f})(s-l^*_D(s), s),
\end{eqnarray*}
where the last equality is obtained by L'H\^{o}pital's rule, and hence $l^*_D(s)$ is the value which gives the supremum to $\frac{\varphi(s)}{\varphi(s-z)}
Q(s;z)
(g-\bar{f})(s-z, s)$.

When $\log \varphi(\cdot)$ is linear,
in lieu of \eqref{eq:V-V-1}, we claim that
\begin{equation}
 \frac{V_\epsilon(s)}{\varphi(s)}=\frac{1}{1+\epsilon}\cdot\frac{V(s+\epsilon)}{\varphi(s+\epsilon)}
\end{equation} for $\epsilon>0$ sufficiently small.
The intuition here is the following: since $\alpha_s(\epsilon)$ in Lemma \ref{lem:epsilon} can be written as  $\left(1+\epsilon\cdot\frac{\varphi''(s)}{\varphi'(s)}\right)\uparrow 1$ (as $\epsilon\downarrow 0$) and in case of Brownian motion, $\frac{\varphi''(s)}{\varphi'(s)}=(\text{const})$, the factor should be independent of $s$. It can be easily seen that the proof of Lemma \ref{lem:epsilon} holds in this case, too.  Accordingly, instead of \eqref{eq:V-epsilon}, we have
\[
V_\epsilon(s)-\frac{\varphi(s)}{\varphi(s+\epsilon)}\exp\left(-\epsilon P(s; l^*_D)\right)V(s+\epsilon, s+\epsilon)=\left(1-\frac{\varphi(s)}{\varphi(s+\epsilon)}(1+\epsilon)\exp\left(-\epsilon P(s; l^*_D)\right)\right)V_\epsilon(s).
\] For the rest, we just proceed as in the proof of Proposition \ref{prop:2} to obtain $\tilde{Q}(\cdot ; l_D)$.
\end{proof}

\begin{remark}\label{rem:meaning}\normalfont  Let us slightly abuse the notation by writing
$\frac{F'(s)\varphi'(s)}{\varphi''(s)[F(s)-F(s-z)]+F'(s)\varphi'(s)}=Q(s; z)$ to avoid the long expression.
Note that {\rm $\frac{\varphi(s)}{\varphi(s-z)}Q(s;z)(g-\bar{f})(s-z,s)$ is the value corresponding to the strategy $D$ with $l_D(s)=z$ and $l_D(m)=l^*_D(m)$ for every $m>s$; that is, this amount is obtained when we stop if $X$ goes below $s-z$ in the excursion at level $S=s$ and behave optimally at all  the higher levels $S>s$.}
\end{remark}

In summary, in Case (1), we should resort to Proposition \ref{prop:2} to compute $V(s, s)$.  If $l^*_D(s)\ge 0$ for this $s$, we shall have $(g-\bar{f})(s, s)\le V(s, s)$.

\subsection{Case (2): \emph{$s\in A_s=[x^*(s), r)$}}
Recall that $A_s$ is defined in Assumption \ref{assumption}.  In this case, similar to Case (1), a positive $l^*_D(s)$ may lead to improvement of the value of $V(s, s)$, so that we use Proposition \ref{prop:2}.

\subsection{Case (3): \emph{$s\in B_s=(l, x^*(s)]$}}
%Denote by $\hat{s}$ the point at which
%\begin{equation}\label{eq:s-hat}
%\hat{s}:= \mbox{argmin} \{ x^*(s): s\in (l, x^*(s)]\}
%\end{equation} holds.  This is the left-end point of the stopping region, that is, $(x^*(\hat{s}), \hat{s})$ is in $\Gamma\subset \R^2$. We clam the following:
%\begin{align}\label{eq:s-x(s)}
%&\text{$\forall s<\hat{s}$,  the point $(x^*(\hat{s}), \hat{s})\in \Gamma\subset \R^2$  is the  point at which a sample path }\\
%&\text{from $(s, s)\in \R^2$ \emph{first} enters the stopping region $\Gamma\subset \R^2$}.\nonumber
%\end{align}
%The reason is the following:  Fix an $s<\hat{s}$ and suppose that $s\le x^*(s)$.  %For concreteness, let us assume that $x^*(s)>\hat{s}$.
%Then the interval $(x^*(\hat{s}), x^*(s))\subset \R_+$ would \emph{not} belong to $\Sigma_s$. On the other hand, $(x^*(\hat{s}), x^*(s))$ is in $\Sigma_{\hat{s}}$.  By the definition of $\hat{s}$, $\Sigma_{\hat{s}}\subset \Gamma$ and thus it must be the case that $(\hat{s}, x^*(s))$ cannot be in the continuation region $\mathrm{C}$: otherwise, we would violate the Markovian character of $(X, S)$.  Hence \eqref{eq:s-x(s)} is true.  One obtains the value of $(g-\bar{f})(x^*(\hat{s})), \hat{s})$ when stops there.

%A typical case is that $\hat{s}$ satisfies
$B_s$ is defined in Assumption \ref{assumption}.  For this case, the typical situation is that $x^*(s)$ is monotonically decreasing in $s$.  See Figure \ref{Fig-example-design}.  The curve separating the region $\Gamma$ and $\mathrm{C}_2$ corresponds to the function $x^*(s)$.  Then define the point  $\hat{s}$  such that
\[
s=x^*(s)
\] holds.  One receives $(g-\bar{f})(x^*(\hat{s}), \hat{s})=(g-\bar{f})(\hat{s}, \hat{s})$ when stops there.    In contrast to the previous Cases (1) and (2), $s$ is located  to the left of $x^*(s)$, % and hence for $u\ge s$, there are no $l^*_D(u)> 0$.  Proposition \ref{prop:2} is  not relevant, so we shall go back to Proposition \ref{prop:1}.  Since
we are not supposed to stop during the excursions from the level $u\in [s, \hat{s})$.  Mathematically, it means that we let $u-l^*_D(u)\rightarrow l$ in \eqref{eq:V(s,s) integral form} of Proposition \ref{prop:1}.  The left boundary $l$ is assumed to be natural and hence $F(u-l_D^*(u))\rightarrow 0$ and
$S_{\tau_{l_D}(m)}\in \diff m=\delta_{\hat{s}}\diff m$, the Dirac measure sitting at $\hat{s}$. Now, from \eqref{eq:reach-m},  \eqref{eq:V(s,s) integral form} simplifies to
\begin{align}\label{eq:V-in-continuation}
  V(s, s)&=\int^\infty_s\frac{\varphi(s)}{\varphi(m)}\exp\left( -\int^m_s\frac{F'(u)}{F(u)}\diff u \right)(g-\bar{f})(m, m)\delta_{\hat{s}}\diff m=\frac{\psi(s)}{\psi(\hat{s})}(g-\bar{f})(\hat{s}, \hat{s}),
\end{align}
which is, in view of \eqref{eq:psi-phi},  simply the expected discounted value of $(g-\bar{f})(\hat{s}, \hat{s})$.  Note that for $s\le \hat{s}$, the reward $(g-\bar{f})(x^*(\hat{s}), \hat{s})$ does not depend on $s$ and thereby
with respect to this reward,  $x^*(s)=x^*(\hat{s})$ for $s\le \hat{s}$.\\ %For the case where there is no $s$ such that \eqref{eq:s-hat}, see Remark \ref{rem:no-fixed-point}. \\

Before moving on to find the general solution $V(x, s)$, it should be beneficial to briefly review some special cases in finding $V(s,s)$.  In this section, the diffusion $X$ is geometric Brownian motion $\diff X_t =\mu X_t\diff t+\sigma X_t\diff B_t$ and $(\A-q)v(x)=0$ provides $\varphi(x)=x^{\gamma_0}$ and $\psi(x)=x^{\gamma_1}$ with $\gamma_0<0$ and $\gamma_1>1$.  The parameters are $(\mu, \sigma, q, K, k)=(0.05, 0.25, 0.15, 5, 0.5)$. The values of the options here are computed under the physical measure $\p$.

\subsubsection{Lookback Option}\label{sec:lookback}
The reward function is $(g-\bar{f})=s-kx$ where $k\in [0, 1]$. Set $s=5$.  By setting $y=F(x)$ in \eqref{eq:H_s}, $H_s(F(x))=\frac{s-kx}{\varphi(x)}$.  The graph of $\frac{s-kx}{\varphi(x)}$ against the horizontal axis $F(x)$ is in Figure \ref{Fig-put}-(a). It can be seen that $s\in \Sigma_s$ and that Case (1) applies.  The optimal threshold $l^*_D(s)$ can be found by Proposition \ref{prop:2}: the optimal level $x^*$ is given by $x^*=\beta s$ where $\beta=0.701636$, independent of $s$, so that $l^*_D(s)=(1-\beta)\cdot s$.

Once $l^*_D(s)$ is obtained, we can compute $V(s, s)$ from \eqref{eq:V-explicit}.  While we shall discuss the general method of computing $V(x, s)$ for $x\le s$ in Section \ref{sec:method}, it is appropriate to touch upon this issue here.  For this fixed $s=5$, we examine the smallest concave majorant of $H_s(y)$.  But the majorant must pass the point \[\left(F(s), \frac{V(s, s)}{\varphi(s)}\right) \conn \Big(F(s-l^*_D(s)), H_s(s-l^*_D(s))\Big).\]  The red line $L_s$ is drawn connecting  these points with a positive slope $\gamma(s-l^*_D(s))$.  In fact, the smooth-fit principle holds at $F(s-l^*_D(s))$ as is discussed in \cite{shepp-shiryaev-1993}.  Accordingly, $(s-l^*_D(s), s)\subset \R$ is in the continuation region.  For more comments about the smooth-fit principle, see Appendix \ref{rem:shiryaev}.

\subsubsection{Perpetual Put} \label{sec:perpetual}
\begin{figure}[h]
\begin{center}
\begin{minipage}{0.45\textwidth}
\centering{\includegraphics[scale=0.75]{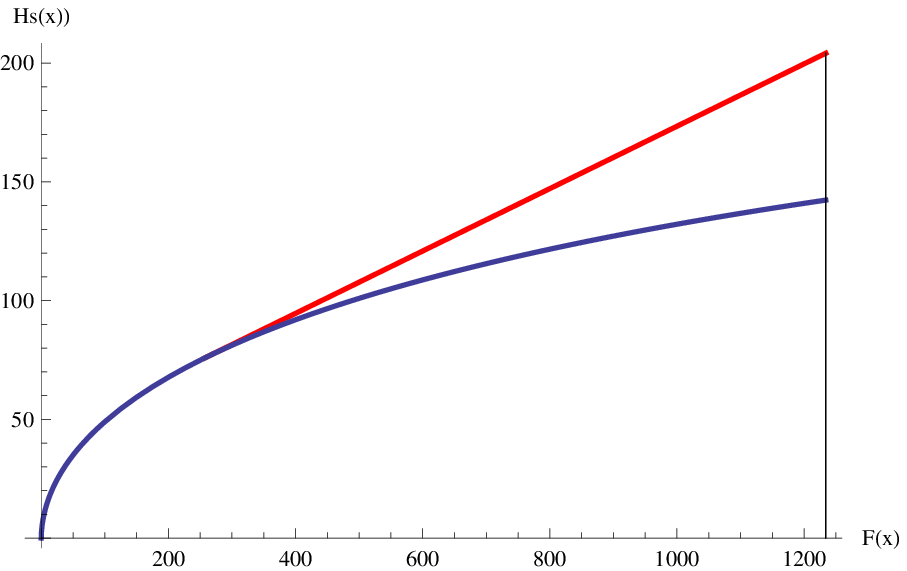}}\\
(a) Lookback Option $(k=1/2)$
\end{minipage}
\begin{minipage}{0.45\textwidth}
\centering{\includegraphics[scale=0.75]{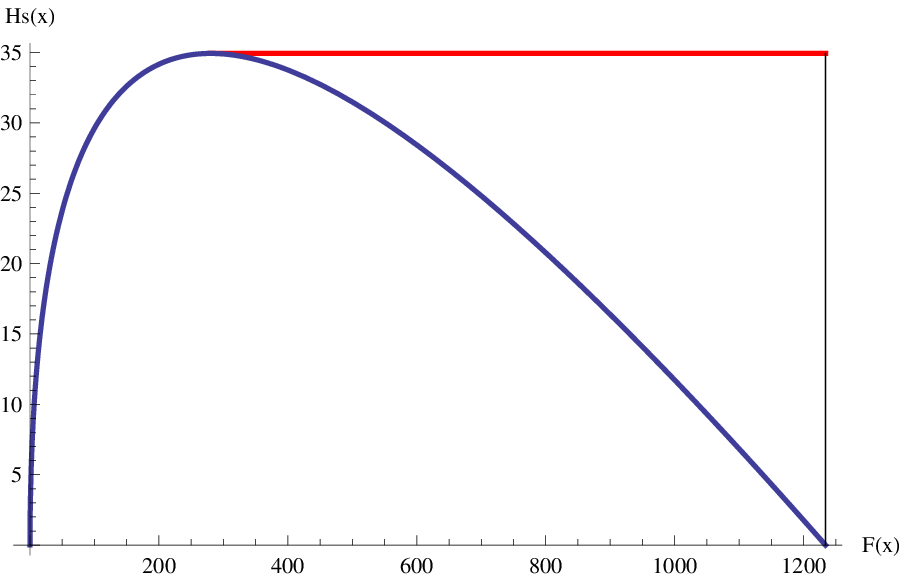}}\\
(b) Perpetual Put $(K=5)$
\end{minipage}
\caption{\small \textbf{The graphs of $g(x, s)/\varphi(x)$ against the horizontal axis $F(x)$ : }  We fix $s=5$ in both problems. The vertical lines show the position of $\varphi(s)$.    In the perpetual put case, the optimal exercise threshold is well-known: $x^*=\frac{\gamma_0K}{\gamma_0-1}=3.57604$, which does not depend on $s$. } \label{Fig-put}
\end{center}
\end{figure}

The reward function is $(g-\bar{f})(x, s)=g(x)=(K-x)^+$ which does not depend on $s$ and there is no absorbing boundary. This is an ordinary optimal stopping problem.   The graph of $g(x)/\varphi(x)$ is in Figure \ref{Fig-put}-(b) which is drawn against the horizontal axis of $y=F(x)$ when $s=5$.  The function
$H_s(F(x))$ attains \emph{unique maximum} at $F(x^*)$ where $x^*=3.57604$, so that $l^*_D(s)=s-x^*=1.42396$.  Since $g$ is independent of $s$, so is $x^*$.

By using this fact, we can use Proposition \ref{prop:1} and compute $V(s, s)$ easily.  In fact, an observation of \eqref{eq:V(s,s) integral form} reveals that if $F(u-l_D(u))$ is constant for $u>s$, then we have
\[\exp\left(-\int_s^\infty \frac{F'(u)}{F(u)-F(u-l_D(u))}\diff u\right)=0,
\] which in turn makes
\[
\int_s^\infty \exp\left( -\int^m_s\frac{F'(u)\diff u}{F(u)-F(u-l_D(u))} \right)\frac{F'(m)}{F(m)-F(m-l_D(m))}\diff m=1.
\] Then \eqref{eq:V(s,s) integral form} reduces to
\begin{equation}\label{eq:simple-case-no-s}
  V(s, s)=\sup_{l_D(s)}\frac{\varphi(s)}{\varphi(s-l_D(s))}(g-\bar{f})(s-l_D(s), s),
\end{equation}
and $l_D^*(s)$ is the maximizer of the map $z\mapsto \frac{\varphi(s)}{\varphi(s-z)}(g-\bar{f})(s-z, s)$.

At this point the tangent line has slope zero; that is, $\gamma(s-l^*_D(s))=0$.  See the red horizontal line connecting two points $(F(x^*), H_s(x))$ and $(F(s), \frac{V(s, s)}{\varphi(s)})$. At $F(x^*)$, we have the smooth-fit principle hold and $(x^*, s)\subset \R$ is in the continuation region.

\section{General Solution}\label{sec:method}
Finally, let us consider the general case, $S_0\geq X_0$. Since we calculated $V(s,s)$, we can represent $V(x,s)$ by (\ref{eq:one-dim-version}):
\begin{eqnarray}\label{eq:gen}
V(x,s)&=&\sup_{\tau\in\S}\E^{x,s}\left[\1_{\{\tau<T_s\}}e^{-q\tau}(g-\bar{f})(X_{\tau},s)
+\1_{\{T_s<\tau\}}e^{-qT_s}V(s,s)\right].
\end{eqnarray}
As we noted in Section 2, this can be seen as just an one-dimensional optimal stopping problem for the process $X$. In terms of the $(s, x)$-diagram like Figure \ref{Fig-example-design}, we fix $s=\bar{s}$, say and use the information of $V(\bar{s}, \bar{s})$, compute $V(x, \bar{s})$ and tell, by moving down from the diagonal point $(\bar{s}, \bar{s})$,   whether a point $(\bar{s}, x)$ belongs to $\mathrm{C}$ or $\Gamma$.   %Since optimal stopping time  is known to be an exit time from a region (Section \ref{sec:diffusion-facts}), we can restrict the set $\S$ of stopping times to $S'$ defined in (\ref{eq:tau-D}).
After discussing generality here, we shall study an example in Section \ref{sec:call-with-s} by showing
how to implement the method presented in Dayanik and Karatzas \cite{DK2003}.

Now suppose that we have found $V(s, s)$ for each $s\in \R_+$.  The next step is to solve \eqref{eq:gen}.  Consider an excursion from the level $S=s$. Recall that $V(s,s)$ represents the value that one would obtain when $X$ would return to that level $s$. If there is no absorbing boundary, we can let the height of excursions arbitrarily large.  Since we are assuming \eqref{eq:finiteness}, by Proposition 5.12 in \cite{DK2003}, the value function in the transformed space must pass the points: \[\Big(0, \xi_l\Big) \conn \left(F(s),\frac{V(s,s)}{\varphi(s)}\right).\]
 Then the task is to find for each $s$ the smallest concave majorant $W_s(y)$ of \[H_s(y):=\frac{(g-\bar{f})(F^{-1}(y),s)}{\varphi(F^{-1}(y))}\] and to identify the region $\{y: H_s(y)=W_s(y)\}$ as optimal stopping region.
Mathematically, $W_s$ must satisfy the following conditions:
\begin{figure}[h]
\begin{center}
\begin{minipage}{0.6\textwidth}
\centering{\includegraphics[scale=1]{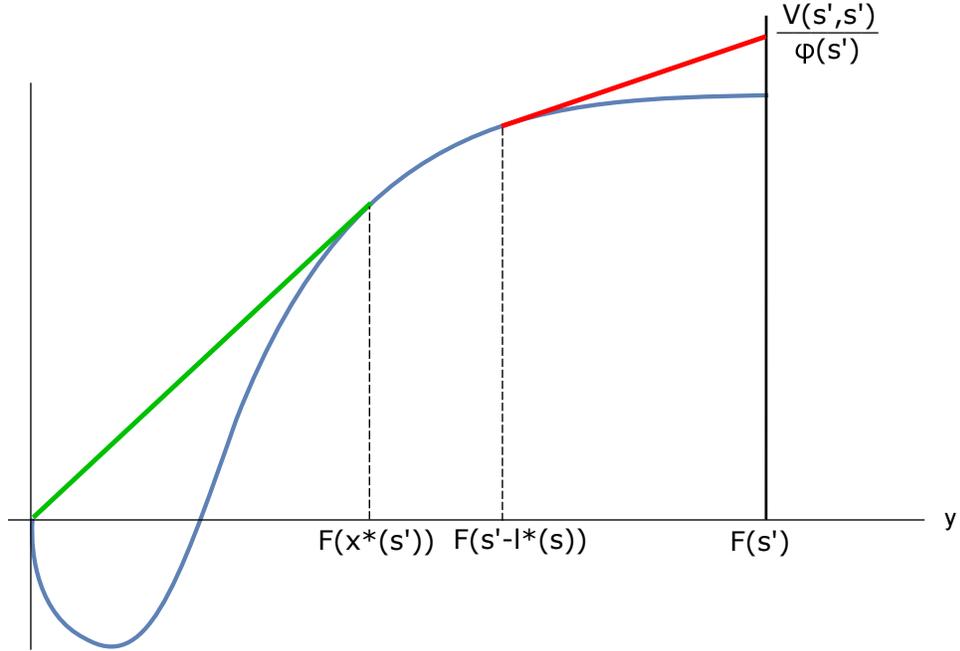}}\\
\end{minipage}
\caption{A typical example of $W_{s}$ and $H_{s}$. \small Fix $s=\s$.  Find $V(s, s)$ and then specify the optimal strategy on $(l, s]$ based on $H_{\s}$.} \label{Fig1}
\end{center}
\end{figure}

\begin{enumerate}
\renewcommand{\labelenumi}{(\roman{enumi})}
 \item $W_s(y)\geq H_s(y)$ on $[0, F(s)]$,
 \item $W_s(F(s))=\frac{V(s,s)}{\varphi(s)}$,
 \item $W_s(0)=\xi_l$,
 \item $W_s$ is concave on $[0, F(s)]$, and
 \item for any functions $\overline{W}_s$ which satisfies four conditions above, $W_s\leq\overline{W}_s$ on $[0, F(s)]$.
\end{enumerate}
Now once we have done with one $s$, we then move on to another $\tilde{s}$, say, and find $W_{\tilde{s}}$ in the new interval
$[0, F(\tilde{s})]$.

Figure \ref{Fig1} illustrates a typical example of the graphs of $W_{s}$ and $H_{s}$ in transformed space.  Fix $s=s'$. Take the point $F(s')$ on the horizontal axis and find $W_{s'}(y)$ that satisfies the above conditions. For this purpose, three vertical lines are drawn at $y=F(x^*(s')), F(s'-l^*_D(s'))$, and $F(s')$ from the left to right. Starting with the point $\left(F(\s), \frac{V(\s, \s)}{\varphi(\s)}\right)$, the concave majorant near that point is the line that is tangent to $H_{s'}(y)$.  The tangency point is $(F(\s-l^*_D(\s), H_{\s}(F(\s-l^*_D(\s)))$.  In the region $[F(x^*(\s), F(\s-l^*_D(\s)]$, the value function is the reward function itself.  On the other hand, the smallest concave majorant of $H_\s(y)$ on $(0, F(x^*(\s))$ is the line, from the origin, tangent to $H_\s(y)$ at $F(x^*(\s))$.

For this $\s$, optimal strategy reads as follows: If it happens that $x \in (0, \s)$ belongs to $(\s-l^*_D(\s), \s)$, one should see if $X$ reaches $\s-l^*_D(\s)$ or $\s$, whichever comes first.  If the former point is the case, one should stop and receive the reward, otherwise one should continue with $s>\s$.  If $x \in (0, \s)$ belongs to $(x^*(\s), \s-l^*_D(\s))$, one should immediately stop $X$ and receive $g(x, \s)$.  Finally, in $x\in (0, x^*(\s))$, one should wait until $X$ reaches $x^*(\s)$.

\subsection{Illustration: A New Problem}\label{sec:call-with-s}
To illustrate how to implement the solution method for a problem that involves both $S$ and $X$, we postulate the reward function as
\begin{equation}\label{eq:example-reward}
g(x, s)=s^a + kx^b-K, \quad a, b, k, K>0
\end{equation} and $f(x, s)\equiv 0$.   For concreteness, we set $a=1/2, b=1, k=1/2$, and $K=5$.
We assume that the underlining process $X$ is  geometric Brownian motion:
\begin{equation*}
\diff X_t=\mu X_t \diff t + \sigma X_t \diff B_t,\quad t\in \R,
\end{equation*}
where $\mu$ and $\sigma$ are constants and $B$ is a standard Brownian motion under $\mathbb{P}$. In this case, \[\psi(x) =x^{\gamma_1} \quad\text{and}\quad \varphi(x)=x^{\gamma_0}
\]
where
\begin{eqnarray*}
\gamma_0&=&\frac{1}{2}\left(-\left(\frac{2\mu}{\sigma^2}-1\right)-\sqrt{\left(\frac{2\mu}{\sigma^2}-1\right)^2+\frac{8q}{\sigma^2}}\right)<0,
\end{eqnarray*}
and
\begin{eqnarray*}
\gamma_1&=&\frac{1}{2}\left(-\left(\frac{2\mu}{\sigma^2}-1\right)+\sqrt{\left(\frac{2\mu}{\sigma^2}-1\right)^2+\frac{8q}{\sigma^2}}\right)>1.
\end{eqnarray*}
\begin{figure}[h]
\begin{center}
\begin{minipage}{0.45\textwidth}
\vspace{0.7cm}
\centering{\includegraphics[scale=0.525]{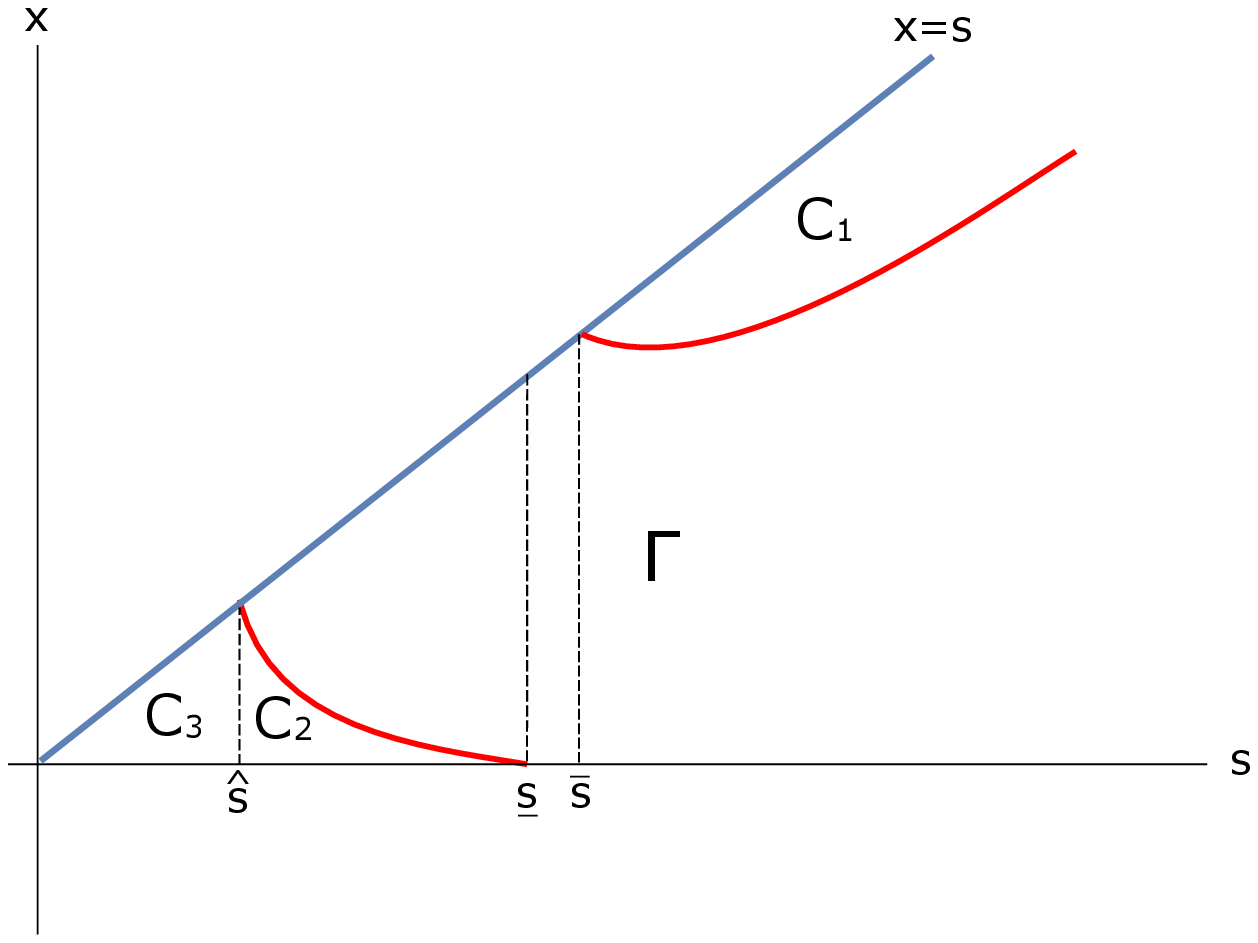}}\\
\end{minipage}
\begin{minipage}{0.45\textwidth}
\centering{\includegraphics[scale=0.55]{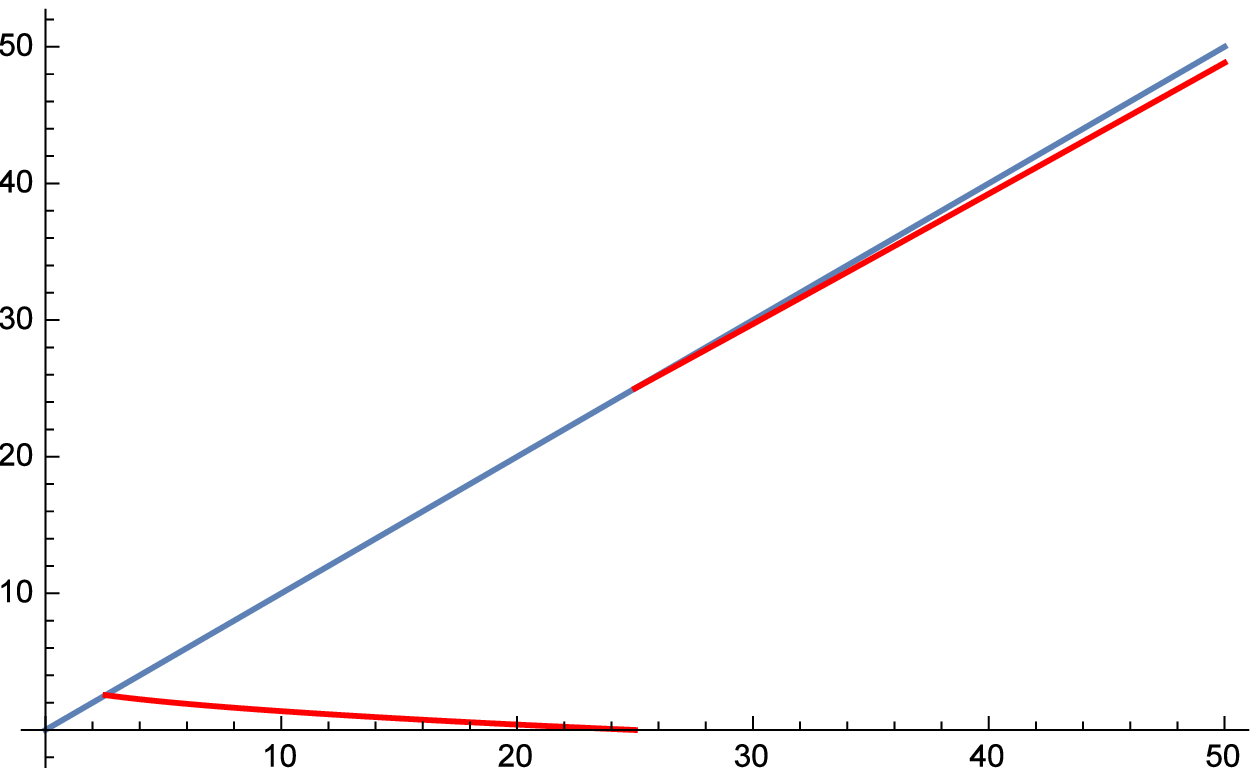}}\\
\end{minipage}
\caption{\small The solution $(s, x)$-diagram to the reward function \eqref{eq:example-reward} in a schematic presentation (left) and in the real scale (right). } \label{Fig-example-design}
\end{center}
\end{figure}
For this reward $g-\bar{f}$, the value of $\xi_l$ in \eqref{eq:finiteness} is zero for any $s\in \II$.

It should be helpful to preview  the entire state space at the beginning.  See Figure \ref{Fig-example-design} for optimal strategy: continuation and stopping region in the two-dimensional diagram.  The left panel is for presenting in a schematic drawing.  For different values of $(s, x)$, we see continuation regions $\mathrm{C}_1, \mathrm{C}_2, \mathrm{C}_3$ and stopping region $\Gamma$. The right panel is the real solution for this problem with parameters $(\mu, \sigma, q)=(0.05, 0.25, 0.15)$.  According to the values of $s$, we have four regions: (i) $s>\bar{s}$,  (iii) $\underline{s}\le s \le \bar{s}$,  (ii) $\hat{s}<s< \underline{s}$,  and (iv) $s\le \hat{s}$. We shall explain how to find the value function and optimal strategy for each region.  Note that for the ease of exposition, we handle $\hat{s}<s< \underline{s}$ before we do $\underline{s}\le s \le \bar{s}$.

\begin{figure}[h]
\begin{center}
\begin{minipage}{0.3\textwidth}
\centering{\includegraphics[scale=0.5]{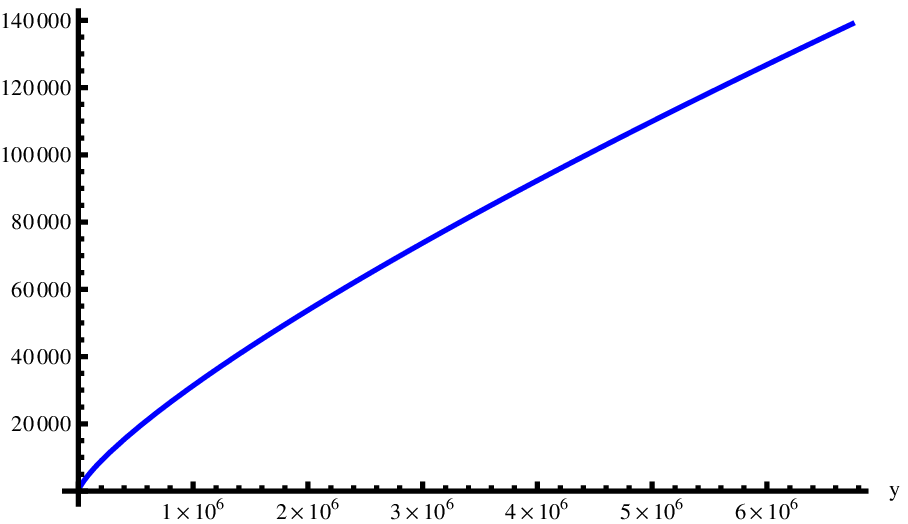}}\\
\small (i) Graph of $H_s(y)$ when $s=35$
\end{minipage}
\begin{minipage}{0.3\textwidth}
\centering{\includegraphics[scale=0.5]{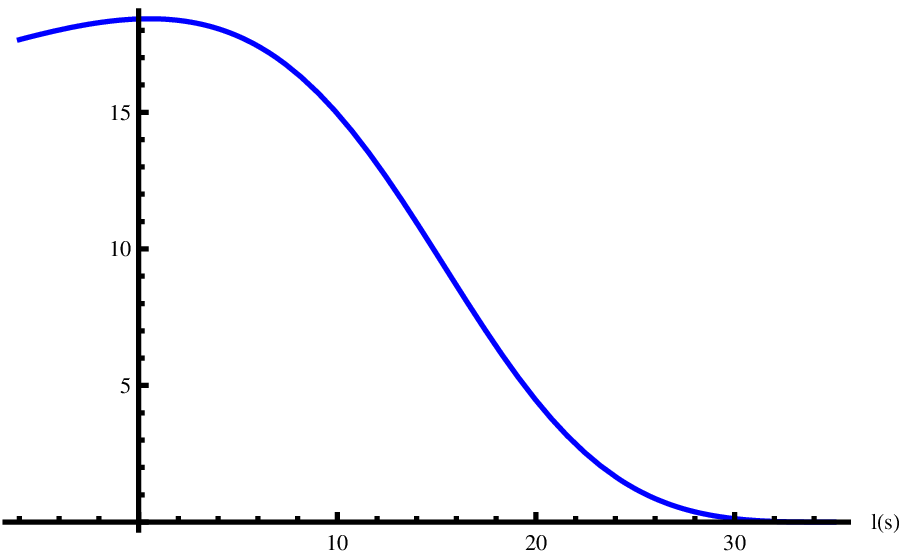}}\\
\small (ii) The map $A(l)$ to find $l^*_D(s)=0.5371.$
\end{minipage}
\begin{minipage}{0.3\textwidth}
\centering{\includegraphics[scale=0.5]{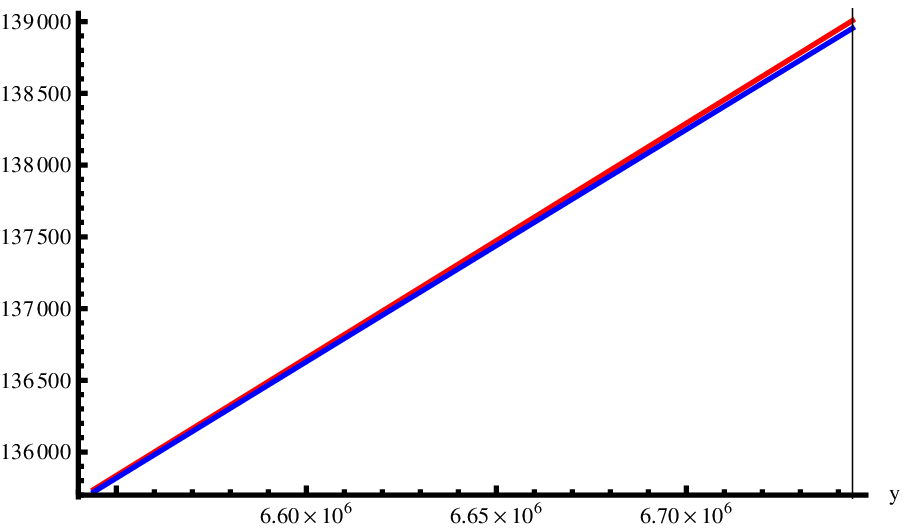}}\\
\small (iii) The value function $V(x, s)$ (red) and $H_s(y)$ (blue)
\end{minipage}
\caption{\small The region $s \in (\bar{s}, +\infty)$. Note in panel (iii) for  a better picture, the lower-left corner of the graph is not the origin.   }
\label{fig:example-region1}
\end{center}
\end{figure}

\noindent (i) Let us start with $s>\bar{s}$.  First we need to find $V(s, s)$. Plug this $s$ in \eqref{eq:example-reward} and examine the reward function in the transformed space:
\[H_s(y)=\frac{(g-\bar{f})(F^{-1}(y), s)}{\varphi(F^{-1}(y))}=\frac{\sqrt{s}+ky^{\frac{1}{\gamma_1-\gamma_0}}-K }{y^{\frac{\gamma_0}{\gamma_1-\gamma_0}}}, \quad y\in (0, \infty).\]
Figure \ref{fig:example-region1} -(i) shows  the function $H_{s}(y)$ on $[0, F(s)]$.  For this $s$, $H_s(y)$ is concave in the neighborhood of $F(s)$ and Case (a) in Section \ref{sec:finding-V(s,s)} applies.  The map
\[A(l):
l \mapsto
\frac{\varphi(s)}{\varphi(s-l)}\cdot\frac{F'(s)\varphi'(s)}{\varphi''(s)[F(s)-F(s-l)]+F'(s)\varphi'(s)}\cdot g(s-l,s)\]
 in \eqref{eq:V-explicit} is in Figure \ref{fig:example-region1}-(ii)
and $l^*_D(s)=0.5371$ when $s=35$.

Now we can find $V(x, s)$ for $x\in [0, s]$.  Following Section \ref{sec:method}, we shall find the smallest concave majorant of $H_{s_1}(y)$ that passes the origin and $\left(F(s),\frac{V(s, s)}{\varphi(s)}\right)$.  For this particular $s$, a diagram similar to Figure \ref{Fig1}-(b) can be drawn: see Figure \ref{fig:example-region1}-(iii).
The value function on the continuation region $\mathrm{C}_1$ (red line in the graph) is
\[
V(x, s)=W(F(x))\varphi(x)= \Big(\beta_1(F(x)-F(s-l^*_D(s)))+H_s(F(s-l^*_D(s)))\Big)\cdot \varphi(x)
\]
where $\beta_1=\frac{\diff H_s(y)}{\diff y}\Bigm|_{y=F(s-l^*_D(s))}$.

In summary, the value function is
\begin{align*}
V(x, s)=\begin{cases}
 g(x, s), &  x\in (0, s-l^*_D(s)),\\
  \Big(\beta_1(F(x)-F(s-l^*_D(s)))+H_s(F(s-l^*_D(s)))\Big)\cdot \varphi(x), & x\in [s-l^*_D(s), s].
\end{cases}
\end{align*}

\begin{figure}[h]
\begin{center}
\begin{minipage}{0.45\textwidth}
\centering{\includegraphics[scale=0.75]{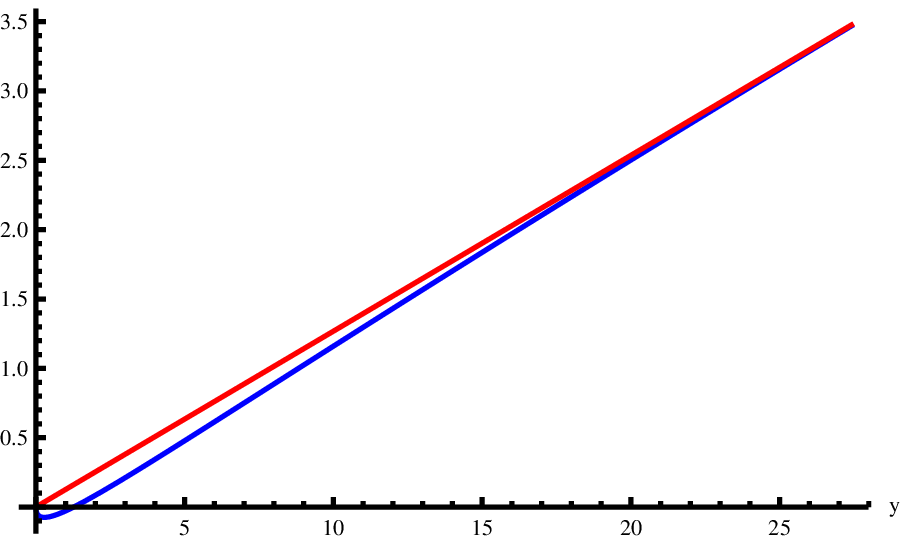}}\\
(i) Graph of $H_s(y)$ and the tangent line $W_s(y)$ when $s=20$. $x^*(s)=2.1242$.
\end{minipage}
\begin{minipage}{0.45\textwidth}
\centering{\includegraphics[scale=0.75]{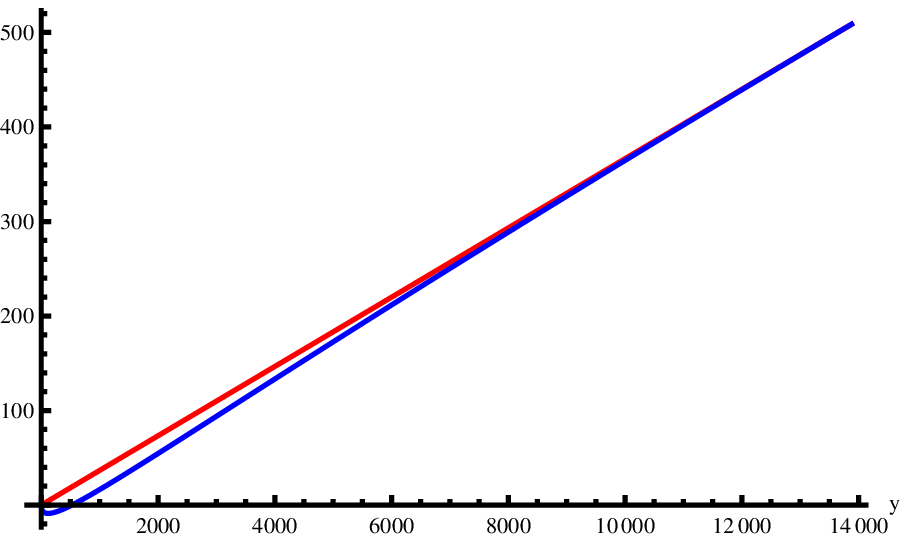}}\\
(ii) Graph of $H_{\hat{s}}(y)$ and the tangent line $W_{\hat{s}}(y)$ with $\hat{s}=x^*(\hat{s})=8.6420$
\end{minipage}
\caption{The region $s \in (\hat{s}, \underline{s})$ (left) and $s\in (0, \bar{s}]$ (right). }
\label{fig:example-region3}
\end{center}
\end{figure}

\noindent (ii) Let us move on to the region $\hat{s}<s< \underline{s}$ (before we examine $\underline{s}\le s \le \bar{s}$).  See Figure \ref{fig:example-region3}-(i) for the graph of $H_s(y)$ in the transformed space.  In the neighborhood of $F(s)$, the reward $H_s(y)$ is concave but We have $l^*_D(s)=0$ for $s$ in this region, so that $V(s, s)=g(s, s)$.  On the other hand, for the reward function $g(x, s)$ with this fixed $s\in (\hat{s}, \underline{s}]$, there is a point $x^*(s)$ such that $(0, x^*(s))$ is the continuation region.  Let us see this situation in the transformed space. See Figure \ref{fig:example-region3}-(i) again.  At the point $F(x^*(s))>F(s)$, the smallest concave majorant $W_s(y)$ of $H_s(y)$ is the line, from the origin, tangent to $H_s(y)$.  Hence the value function is
\begin{align*}
V(x, s)=\begin{cases}
  (\beta_2F(x))\cdot \varphi(x)=\beta_2\psi(x) , &  x\in (0, x^*(s)],\\
g(x, s), & x\in (x^*(s), s],
\end{cases}
\end{align*}
where $\beta_2=\frac{\diff H_s(y)}{\diff y}\Bigm|_{y=F(x^*(s))}$. Hence $\beta_2\psi(x)$ is the value function in region $\mathrm{C}_2$ in Figure \ref{Fig-example-design}.

\noindent (iii) For $\underline{s}\le s \le \bar{s}$, we have $l^*_D(s)=0$ so that the point $(s, s)$ is in the stopping region.  Moreover, there exist no points $x^*(s)$ where the line from the origin becomes tangent to $H_s(y)$.  Accordingly, the value function is
\[
V(x, s)=g(x, s)\quad x \in (0, s].
\]
In our parameters, $\underline{s}=25$ and $\bar{s}$ is very close; $\bar{s}>25$.

\noindent (iv) Now we shall examine $s\le \hat{s}$.  As explained in Case (3) in Section \ref{sec:finding-V(s,s)}, this $\hat{s}=8.6420$ satisfies $\hat{s}=x^*(\hat{s})$.  See \eqref{eq:V-in-continuation}.  Following the argument there, $V(s, s)=\frac{\psi(s)}{\psi(\hat{s})}(g-\bar{f})(x^*(\hat{s}), \hat{s})$.  Accordingly, the value function is
\begin{equation}\label{eq:final}
V(x, s)=\frac{\psi(x)}{\psi(\hat{s})}(g-\bar{f})(x^*(\hat{s}), \hat{s}), \quad x\in (0, s].
\end{equation}
This corresponds to region $\mathrm{C}_3$ in Figure \ref{Fig-example-design}.

%\begin{remark}\normalfont\label{rem:no-fixed-point}
%Suppose that there is no $s$ such that $s=x^*(s)$ holds.  Then we can instead define
%\[
%\hat{s}:=\inf\{s: (s, x)\in \Gamma \quad \text{for some}\word x<s\}.
%\] Then we have
%\[
%V(x, s)=\frac{\psi(x)}{\psi(\hat{s})}(g-\bar{f})(x^*(\hat{s}), \hat{s}),
%\] which is essentially the same as \eqref{eq:final}.  Recall that $x^*(\cdot)$ is defined in Section \ref{sec:finding-V(s,s)}.
%\end{remark}

\begin{appendix}
  \section{Proof of Lemma \ref{lem:epsilon}} \label{app:A}
Lemma \ref{lem:epsilon} is the following claim: \emph{Under the assumption of Proposition \ref{prop:2} with convex $\log \varphi(\cdot)$, we have
\begin{equation} \label{eq:V-V}
    \frac{V_\epsilon(s)}{\varphi(s)}=\alpha_s(\epsilon)\cdot\frac{V(s+\epsilon, s+\epsilon)}{\varphi(s+\epsilon)} \quad \word\text{where}\word
    \alpha_s(\epsilon):=\frac{\varphi'(s+\epsilon)}{\varphi'(s)}.
\end{equation}}
\begin{proof}(of the lemma)
Recall \eqref{eq:lD} for the definition of $l_D(s)$.
In view of \eqref{eq:probability}, the probabilistic meaning of \eqref{eq:V-epsilon} is that $V_\epsilon(s)$ is attained when one chooses the excursion level $l_D(s)$ optimally in the following optimal stopping:
\begin{equation}\label{eq:new-osp}
  V_\epsilon(s)=\sup_{l_D(s)}\E^{s, s}[e^{-qT_{s+\epsilon}}\1_{\{T_{s+\epsilon}\le \tau_{s-l_D(s)}\}}V(s+\epsilon, s+\epsilon)+e^{-q\tau_{s-l_D(s)}}\1_{\{T_{s+\epsilon}> \tau_{s-l_D(s)}\}}(g-\bar{f})(s-l_D(s), s)],
\end{equation}
that is, if the excursion from $s$ does not reach the level of $l_D(s)$ before $X$ reaches $s+\epsilon$, one shall receive $V(s+\epsilon, s+\epsilon)$ and otherwise, one shall receive the reward. By using the transformation \eqref{eq:F-trans}, one needs to consider the function $\frac{(g-\bar{f})(x, s)}{\varphi(x)}$ and the point $\left(F(s+\epsilon), \frac{V(s+\epsilon, s+\epsilon)}{\varphi(s+\epsilon)}\right)$ in the $(F(x), z(x)/\varphi(x))$-plane. Then the value function of \eqref{eq:new-osp} in this plane is the smallest concave majorant of $\frac{(g-\bar{f})(x, s)}{\varphi(x)}$ which passes through the point $\left(F(s+\epsilon), \frac{V(s+\epsilon, s+\epsilon)}{\varphi(s+\epsilon)}\right)$. It follows that $\frac{V_\epsilon(s)}{\varphi(s)}\le  \frac{V(s+\epsilon, s+\epsilon)}{\varphi(s+\epsilon)}$.  As $\epsilon\downarrow 0$, it is clear that $\frac{V(s+\epsilon, s+\epsilon)}{\varphi(s+\epsilon)}\downarrow \frac{V(s, s)}{\varphi(s)}$ and $\frac{V_\epsilon(s)}{\varphi(s)}\downarrow \frac{V(s,s)}{\varphi(s)}$. Suppose, for  a contradiction, that we  have
\begin{equation}\label{eq:for-contradiction}
\alpha_s(\epsilon)\frac{V(s+\epsilon, s+\epsilon)}{\varphi(s+\epsilon)}<\frac{V_\epsilon(s)}{\varphi(s)}<\frac{V(s+\epsilon, s+\epsilon)}{\varphi(s+\epsilon)},
\end{equation} for all $\epsilon>0$.
This implies that the first term in \eqref{eq:for-contradiction} goes to $\frac{V(s,s)}{\varphi(s)}$ from below and the third term goes to the same limit from above.  While the second inequality always hold, the first inequality leads to a contradiction to the fact that the function $\epsilon\mapsto (1-\alpha_s(\epsilon))\frac{V(s+\epsilon, s+\epsilon)}{\varphi(s+\epsilon)}$ is continuous for all $s$.

Indeed, due to the monotonicity of $\alpha_s(\epsilon)\frac{V(s+\epsilon, s+\epsilon)}{\varphi(s+\epsilon)}$ in $\epsilon$, we would have $\frac{V_\epsilon(s)}{\varphi(s)}>\frac{V(s, s)}{\varphi(s)}>\alpha_s(\epsilon)\frac{V(s+\epsilon, s+\epsilon)}{\varphi(s+\epsilon)}$ for all $\epsilon>0$.  Hence one cannot make the distance between $\frac{V(s+\epsilon, s+\epsilon)}{\varphi(s+\epsilon)}$ and $\alpha_s(\epsilon)\frac{V(s+\epsilon, s+\epsilon)}{\varphi(s+\epsilon)}$ arbitrarily small without violating \eqref{eq:for-contradiction}.
This shows that there exists an $\epsilon'=\epsilon'(s)$ such that $\epsilon<\epsilon'$ implies that $\frac{V_\epsilon(s)}{\varphi(s)}\le \alpha_s(\epsilon)\frac{V(s+\epsilon, s+\epsilon)}{\varphi(s+\epsilon)}$.

On the other hand,  in \eqref{eq:new-osp}, one could choose a stopping time $\tau_{l_{D}(s)}$ that visits the left boundary $l$, then by reading \eqref{eq:probability} with $l_D(u)=u$ and $m=s+\epsilon$,
 \eqref{eq:new-osp} becomes
\begin{align*}
V_\epsilon(s)&\ge\frac{\varphi(s)}{\varphi(s+\epsilon)}\exp\left(-\int_s^{s+\epsilon}\frac{F'(u)\diff u}{F(u)-F(0+)}\right)V(s+\epsilon, s+\epsilon)\\
&>\frac{\varphi(s)}{\varphi(s+\epsilon)}\exp\left(-\int_s^{s+\epsilon}\frac{F'(u)\diff u}{F(u)}\right)V(s+\epsilon, s+\epsilon)\\
 &>\frac{\varphi(s)}{\varphi(s+\epsilon)}\frac{F(s)}{F(s+\epsilon)}V(s+\epsilon, s+\epsilon)=\frac{\psi(s)}{\psi(s+\epsilon)}V(s+\epsilon, s+\epsilon)=\E^{s, s}(e^{-qT_{s+\epsilon}})V(s+\epsilon, s+\epsilon)
\end{align*} for \emph{any} $\epsilon>0$.
Since $s\in \mathcal{I}$ is a regular point, the last expectation can be arbitrarily close to unity, monotonically in $\epsilon$ (see page 89 \cite{IM1974}).  Now suppose that there were no $\epsilon$'s such that $V_\epsilon(s)\ge V(s+\epsilon, s+\epsilon)$.  It follows that for any $\epsilon$, we would have
\[
V_\epsilon(s)>\E^{s, s}(e^{-qT_{s+\epsilon}})V(s+\epsilon, s+\epsilon)>\E^{s, s}(e^{-qT_{s+\epsilon}})V_\epsilon(s).
\] Then by letting $\epsilon\downarrow 0$, it would be $V(s)>V(s)$ for all $s\in \II$, which is absurd.  Since the convergence of $\E^{s, s}(e^{-qT_{s+\epsilon}})\uparrow 1$ is monotone in $\epsilon$, there exists an $\epsilon''=\epsilon''(s)>0$ such that $\epsilon<\epsilon''$ implies that $V_\epsilon(s)\ge V(s+\epsilon, s+\epsilon)$.  By using the second assumption in the statement of Proposition, in particular \eqref{eq:log-varphi}, for any $s$, we have
$\frac{V_\epsilon(s)}{\varphi(s)}\ge \alpha_s(\epsilon)\frac{V(s+\epsilon, s+\epsilon)}{\varphi(s+\epsilon)}$ for $\epsilon<\epsilon''$. This completes the proof of Lemma \ref{lem:epsilon}.
\end{proof}
\section{On the smooth-fit principle}\label{rem:shiryaev}
Reference is made to Section \ref{sec:finding-V(s,s)}.  Another representation of $V(s, s)$ is still possible. Continue to fix $s\in \II$.  The smooth-fit principle is assumed to hold at an optimal
point $s-l^*_D(s)$ and we have $(s-l_D^*(s), s)$ as a continuation region with the line $L_s$ tangent to the function $(g-\bar{f})(s-l_D(s), s)/\varphi(s-l^*_D(s))$ at $y=F(s-l^*_D(s))$.  Then we have a first-order approximation of $\frac{V(s+\epsilon, s+\epsilon)}{\varphi(s+\epsilon)}$:
    \begin{equation}\label{eq:approx}
    \frac{V(s+\epsilon, s+\epsilon)}{\varphi(s+\epsilon)}=\frac{V(s, s)}{\varphi(s)}+\gamma(s-l^*_D(s))\cdot \epsilon F'(s)
    \end{equation} where \begin{equation*}\label{eq:gamma}
    \gamma(s-l^*_D(s)):=\frac{1}{F'(s-l^*_D(s))}\left(\frac{(g-\bar{f})(s-l^*_D(s), s)}{\varphi(s-l^*_D(s))}\right)'\ge 0,
    \end{equation*}
    the slope of $L_s(y)$ at $F(s-l^*_D(s))$. See \eqref{eq:derivative}.
Equate equation \eqref{eq:approx} with \eqref{eq:V-epsilon-in-F} to obtain, for any $\epsilon$ sufficiently close to zero,
\[
\left(\frac{V_\epsilon(s)}{\varphi(s)}-\frac{\epsilon P(s; l^*_D)(g-\bar{f})(s- l^*_D(s), s)}{\varphi(s-l^*_D(s))}\right)e^{\epsilon P(s; l^*_D)}
=\frac{V(s, s)}{\varphi(s)}+\gamma(s-l^*_D(s))\cdot \epsilon F'(s).
\] Since $V_\epsilon(s)\mapsto V(s)$ as $\epsilon\rightarrow 0$,
\begin{align}\label{eq:for-special-V}
  \frac{V(s, s)}{\varphi(s)}&=\lim_{\epsilon\downarrow 0} \frac{\gamma(s-l^*_D(s))\cdot \epsilon F'(s)+\epsilon e^{\epsilon P(s; l^*_D)} P(s; l^*_D)\frac{(g-\bar{f})(s- l^*_D(s), s)}{\varphi(s-l^*_D(s))} }{e^{\epsilon P(s; l^*_D)}-1}\nonumber\\
  &=\frac{\gamma(s-l^*_D(s)) F'(s)+ P(s; l^*_D)\frac{(g-\bar{f})(s- l^*_D(s), s)}{\varphi(s-l^*_D(s))} }{P(s; l^*_D)}
  %&=\frac{\varphi(s)}{P(s; l^*_D(s))}\left[\frac{F'(s)}{F'(s-l^*(D))}\left(\frac{(g-\bar{f})(s-l_D^*(s), s)}{\varphi(s-l^*_D(s))}\right)'+P(s; l^*_D(s))\frac{(g-\bar{f})(s-l_D^*(s), s)}{\varphi(s-l_D^*(s))}\right]
\end{align}
where the last equality is obtained by L'H\^{o}pital's rule. By multiplying $\varphi(s)$ on both sides and performing some algebra, we have
\begin{align}\label{eq:V-explicit-2}
V(s,s)&=\frac{\varphi(s)}{P(s; l^*_D(s))}\left[\frac{F'(s)}{F'(s-l^*(D))}\left(\frac{(g-\bar{f})(s-l_D^*(s), s)}{\varphi(s-l^*_D(s))}\right)'+P(s; l^*_D(s))\frac{(g-\bar{f})(s-l_D^*(s), s)}{\varphi(s-l_D^*(s))}\right].
\end{align}

Let us mention some comments on the examples in Section \ref{sec:finding-V(s,s)}. In Section \ref{sec:perpetual} of the perpetual put, for the $l^*_D(s)$, as can be seen from the Figure \ref{Fig-put}, $V(s+\epsilon)/\varphi(s+\epsilon)=V(s, s)/\varphi(s)$.  Then   \eqref{eq:V-explicit-2} (see also \eqref{eq:for-special-V}) reduces to a very simple form
\begin{equation*}
V(s, s)=\frac{\varphi(s)}{\varphi(s-l^*_D(s))}(g-\bar{f})(s-l^*_D(s), s),
\end{equation*}
which is the same as \eqref{eq:simple-case-no-s}.

In the lookback option in Section \ref{sec:lookback}, the case of $k=0$ was solved in the well-known \cite{shepp-shiryaev-1993}. For $k=0$, the reward $(g-\bar{f})(x, s)$ does not depend on $x$, so that we have further simplification of $V(s,s)$.  A straightforward computation from \eqref{eq:V-explicit-2} below %(or \eqref{eq:V-explicit2})
yields
\begin{align}\label{eq:simple-case-no-x}
  V(s,s)=\frac{\psi'(s-l^*)\varphi(s)-\psi(s)\varphi'(s-l^*)}{\psi'(s-l^*)\varphi(s-l^*)-\psi(s-l^*)\varphi'(s-l^*)}\cdot (g-\bar{f})(s)
\end{align}
where we write $l^*:=l^*_D(s)$ for simplicity.  It can be confirmed that \eqref{eq:simple-case-no-x} for $(g-\bar{f})(s)=s$ is the same as
\[
V(s, s)=\frac{s}{\gamma_1-\gamma_0}\left(\gamma_1\left(\frac{1}{\beta}\right)^{\gamma_0}-\gamma_0\left(\frac{1}{\beta}\right)^{\gamma_1}\right)
\] in \cite{shepp-shiryaev-1993}.

We comment on the assumption of the smooth-fit.  When the reward function $(g-\bar{f})(x, s)$ is a \emph{nondecreasing function of the second argument and is differentiable in both arguments},  we have $V(s+\epsilon, s+\epsilon)\ge V(s, s)$ and hence the approximation argument \eqref{eq:approx} is valid.  Note that in Figure \ref{Fig-put}, we have $\frac{V(s, s)}{\varphi(s)}\ge \frac{(g-\bar{f})(s-l^*_D(s), s)}{\varphi(s-l^*_D(s))}$  and the smooth-fit holds at $s-l^*_D(s)$ in both cases.

\end{appendix}

\bibliographystyle{plain}
{\small \bibliography{ospbib2}}

\def\cprime{$'$}
\begin{thebibliography}{10}

\bibitem{alvarez2}
L.~H.~R. Alvarez.
\newblock On the properties of r-excessive mappings for a class of diffusions.
\newblock {\em Ann. Appl. Probab.}, 13 (4):1517--1533, 2003.

\bibitem{alvarez-matomaki2014}
L.~H.~R. Alvarez and P.~Matom\"{a}ki.
\newblock Optimal stopping of the maximum process.
\newblock {\em J. Appl. Prob.}, 51:818--836, 2014.

\bibitem{avram-et-al-2004}
F.~Avram, A.~E. Kyprianou, and M.~R. Pistorius.
\newblock Exit problems for spectrally negative {L}\'{e}vy processes and
  applications to ({C}anadized) {R}ussion options.
\newblock {\em Ann. Appl. Probab.}, 14:215--235, 2004.

\bibitem{Bertoin_1996}
J.~Bertoin.
\newblock {\em {L}\'evy processes}, volume 121 of {\em Cambridge Tracts in
  Mathematics}.
\newblock Cambridge University Press, Cambridge, 1996.

\bibitem{borodina-salminen}
A.~N. Borodin and P.~Salminen.
\newblock {\em Handbook of Brownian Motion - Facts and Formulae, Second
  Edition}.
\newblock Birkh\"{a}user, Basel, Boston, Berlin, 2002.

\bibitem{cinlar}
E.~\c{C}inlar.
\newblock {\em Probability and Stochastics}, volume 261 of {\em Graduate Texts
  in Mathematics}.
\newblock Springer, Berlin Heidelberg, 2011.

\bibitem{DK2003}
S.~Dayanik and I.~Karatzas.
\newblock On the optimal stopping problem for one-dimensional diffusions.
\newblock {\em Stochastic Process. Appl.}, 107 (2):173--212, 2003.

\bibitem{Doney_2005}
R.~A. Doney.
\newblock Some excursion calculations for spectrally one-sided {L}\'evy
  processes.
\newblock In {\em S\'eminaire de {P}robabilit\'es {XXXVIII}}, volume 1857 of
  {\em Lecture Notes in Math.}, pages 5--15. Springer, Berlin, 2005.

\bibitem{dynkin}
E.~B. Dynkin.
\newblock {\em Markov Processes II}.
\newblock Springer, Berlin Heidelberg, 1965.

\bibitem{Guo-Zervos_2010}
X.~Guo and M.~Zervos.
\newblock $\pi$ options.
\newblock {\em Stoch. Proc. Appl.}, 120:1033--1059, 2007.

\bibitem{IM1974}
K.~It\^{o} and H.~P.~McKean Jr.
\newblock {\em Diffusion Processes and their Sample Paths}.
\newblock Springer, Berlin Heidelberg, 1974.

\bibitem{lebedev}
N.~N. Lebedev.
\newblock {\em Special Functions and their Applications, revised edn}.
\newblock Dover Publications, New York, 1972.

\bibitem{ott_2013}
C.~Ott.
\newblock Optimal stopping problems for the maximum process with upper and
  lower caps.
\newblock {\em Ann. Appl. Probab.}, 23:2327--2356, 2013.

\bibitem{peskir1998}
G.~Peskir.
\newblock Optimal stopping of the maximum process: the maximality principle.
\newblock {\em Ann. Probab.}, 26:1614--1640, 1998.

\bibitem{Pham-book}
H.~Pham.
\newblock {\em Continuous-time Stochastic Control and Optimization with
  Financial Applications}, volume~61 of {\em Stochastic Modeling and Applied
  Probability}.
\newblock Springer, Berlin Heidelberg, 2009.

\bibitem{Pistorius_2004}
M.~R. Pistorius.
\newblock On exit and ergodicity of the spectrally one-sided {L}\'evy process
  reflected at its infimum.
\newblock {\em J. Theoretical Probab.}, 17:183--220, 2004.

\bibitem{Pistorius_2007}
M.~R. Pistorius.
\newblock An excursion-theoretical approach to some boundary crossing problems
  and the {S}korokhod embedding for reflected {L}\'evy processes.
\newblock In {\em S\'eminaire de {P}robabilit\'es {XL}}, volume 1899 of {\em
  Lecture Notes in Math.}, pages 278--307. Springer, Berlin, 2007.

\bibitem{shepp-shiryaev-1993}
L.~Shepp and A.~N. Shiryaev.
\newblock The {R}ussian otpion: reduced regret.
\newblock {\em Ann. Appl. Probab.}, 3:631--640, 1993.

\end{thebibliography}

\end{document}